\documentclass[a4paper,11pt]{amsart}

\usepackage{ psfrag, epsfig,amsmath,amssymb,amsthm,array,multicol,verbatim,moreverb,mathrsfs}

\newtheorem{lem}{Lemma}[section]
\numberwithin{equation}{section} \theoremstyle{plain}
\newtheorem{theorem}{Theorem}[section]

\newtheorem{proposition}{Proposition}[section]
\newtheorem{prop}{Proposition}[section]

\newtheorem{rem}{Remark}[section]
\newenvironment{remark}{\begin{rem}\rm}{\smallskip\end{rem}}
\newtheorem*{ex}{Example}

\def\esp{\mathbb{E}}
\def\prob{\mathbb{P}}

\def\real{\mathbb{R}}
\def\Dbb{\mathbb{D}}
\def\Dsf{\mathsf{D}}
\def\Isf{\mathsf{I}}
\def\Hsf{\mathsf{H}}
\def\Cbb{\mathbb{C}}
\def\Mbb{\mathbb{M}}

\def\nat{\mathbb{N}}

\def\Ncal{\mathcal{N}}

\def\Scal{\mathcal{S}}

\def\Xcal{\mathcal{X}}

\def\Dcal{\mathcal{D}}

\def\Ecal{\mathcal{E}}
\def\Fcal{\mathcal{F}}

\def\Rcal{\mathcal{R}}

\def\Cbf{{\bf C}}
\def\Ybf{{\bf Y}}
\def\Xbf{{\bf X}}
\def\Nbf{{\bf N}}
\def\pbf{{\bf p}}
\def\nbf{{\bf n}}
\def\x{{\bf x}}
\def\xbf{{\bf x}}

\def\D{{\bf D}}
\def\A{{\bf A}}
\def\Q{{\bf Q}}
\def\JN{{\bf Net}}
\def\Pbf{{\bf P}}
\def\Xbf{{\bf X}}
\def\Lbf{{\bf L}}

\def\Sbf{{\bf S}}
\def\sbf{{\bf s}}
\def\Pbf{{\bf P}}
\def\ind{{\rm 1\hspace{-0.90ex}1}}

\def\1{{\bf 1}}
\def\x{{\bf x}}

\def\0{{\bf 0}}
\def \essinf{\mathop{\rm essinf}}
\def \esssup{\mathop{\rm esssup}}

\def \Tcal{{\mathcal T}}

\def\Fcal{\mathcal{F}}

\author[Marc Lelarge]{Marc Lelarge}
\title[Large deviations for queueing networks]{Sample path large deviations for queueing networks with Bernoulli routing}
\date{\today}

\begin{document}
\maketitle
\begin{abstract}
This paper is devoted to the problem of sample path large deviations
for multidimensional queueing models with feedback. We derive a new
version of the contraction principle where the continuous map is not
well-defined on the whole space: we give conditions under which it
allows to identify the rate function. We illustrate our technique by
deriving a large deviation principle for a class of networks that
contains the classical Jackson networks.
\end{abstract}
\setcounter{tocdepth}{1} \tableofcontents

%%%%%%%%%%%%%%%%%%%%%%%%%%%%%%%%%%%%%%%%%%%%%%%%%%%%%%%%%%%%%%%%%%%%%%%%%%%%%%%%%%%%%%%%%%%%%%%%%%%%%%%%%%%%%%%%%%
%%%%%%%%%%%%%%%%%%%%%%%%%%%%%%%%%%%%%%%%%%%%%%%%%%%%%%%%%%%%%%%%%%%%%%%%%%%%%%%%%%%%%%%%%%%%%%%%%%%%%%%%%%%%%%%%%%
\newpage

\section*{Introduction}

This paper is concerned with the theory of large deviations of
stochastic processes related to discrete event systems. As opposed
to classical stochastic dynamical systems, for which the evolution
is continuous and described by a stochastic differential equation,
discrete event systems are characterized by synchronization
mechanisms that prevent most of the classical tools to apply. We
present here a new approach for the analysis of the sample path
large deviations of such processes. Unlike standard methods that
require establishing upper and lower bounds, our method relies on
the analogy between the theory of weak convergence and the theory of
large deviations. This analogy is well-known and has been studied by
many authors, we refer to the recent book of Feng and Kurtz
\cite{fenkur} that surveys this field. We should in particular quote
the work of Puhalskii \cite{Puh:pro}, \cite{Puh:book} quite similar
to our approach. However, we will not use the framework of
idempotent measures developed by Puhalskii. We discuss in more
details our general methodology and its relation with the existing
literature after the description of the queueing networks we
consider.

To apply our method we choose a class of queueing networks with
Bernoulli routing, where feedback is allowed. The discontinuous
dynamic of queueing networks makes it hard to study and large
deviations results in the literature are treated on a case by case
basis as in the work of Ganesh and Anantharam \cite{ganan},
Bertsimas, Paschalidis and Tsitsiklis \cite{bpt} or Ramanan and
Dupuis \cite{ramdup}. As we will see, adding the possibility of
feedback makes the problem much harder. For queueing networks with
feedback, existing large deviations results are restricted to
networks described by finite-dimensional Markov processes, see the
works of Dupuis, Ellis and Weiss \cite{dupelwei}, Dupuis and Ellis
\cite{dupel} and Ignatiouk-Robert \cite{ignaop}. In this paper, we
consider networks where the output process of a queue is modeled by
a reflection mapping. This class contains the classical Jackson
networks and our large deviations results extend existing results
for this class obtained by Atar and Dupuis \cite{atardupuis} and
Ignatiouk-Robert \cite{ign}. Our technique allows to obtain large
deviations results under non-exponential assumptions. This case
corresponds to networks with autonomous service and gives an
approximation for queueing networks where each station acts as a
standard single server queue. While preparing this paper, the author
became aware of the work of Puhalskii \cite{Puh:jack} who considers
generalized Jackson networks. The form of the rate function for the
queue length process obtained in \cite{Puh:jack} coincides with our
result, which confirms the intuition that in the large deviation
regime, networks with autonomous service approximate well
generalized Jackson networks. We will discuss more carefully this
result in Section \ref{sec:result}.

In the next section, we give an overview of the general methodology
and then introduce the general notation. Section \ref{sec:contr}
gives an extension of the contraction principle that will allow us
to identify the rate function. Our result is stated without any
reference to any specific discrete event system and could be applied
to other systems. In Sections \ref{sec:extpsi} and \ref{sec:ldp}, we
apply our method to the case of a queueing network.

\subsection*{General methodology}

For simplicity, we adopt here the notation corresponding to our
example of queueing network. As in \cite{lel} or \cite{maj-fluid},
we define the arrival and departure processes $\A$ and $\D$ of each
station of the network as the solution of the fixed point equation
\begin{eqnarray}
\label{fixpointeq}\left\{\begin{array}{lcl}\A &=& \Gamma(\D,\JN),\\
\D &=& \Phi(\A,\JN),\end{array}\right.&\Leftrightarrow &
(\A,\D)=:\Psi(\JN).
\end{eqnarray}
Here $\JN$ is a process that describes all the primitives of the
network as service times at the different stations, routing
decisions, arrival times in the network. The maps $\Gamma$ and
$\Phi$ describe the dynamic of the network.

We consider a sequence of queueing networks $\{\JN_n\}_n$, where the
primitives are counting processes (i.e. $\JN_n$ belongs to a space
denoted by $\Ecal$) and satisfy a large deviation principle (LDP).
We denote by $\Isf^\JN$ the corresponding rate function. It is known
that the map $\Psi$ is well defined if the primitives of the network
are counting processes, see \cite{cm:disc} or \cite{lel} and we
denote $(\A_n,\D_n)=\Psi(\JN_n)$. It is natural to ask whether
$\Psi$ is well defined for processes in $\Dbb$, the space of cadlag
non-decresing functions or at least for absolutely continuous
functions. If this was true, and if $\Psi$ was shown to be
continuous, then we would get thanks to the contraction principle
that the sequence of processes $\{(\A_n,\D_n)\}_n$ satisfies a LDP
with good rate function
\begin{eqnarray}
\label{eq:quote}"\Isf^{\A,\D}(\A,\D)=\inf\left\{ \Isf^{\JN}(\JN),\:
\Psi(\JN)=(\A,\D)\right\}."
\end{eqnarray}
However, the map $\Psi$ turns out not to be well defined for all
possible limits of a sequence of networks $\{\JN_n\}_n\in
\Ecal^\nat$ as defined previously. In particular, the fixed point
equation (\ref{fixpointeq}) can very well be stated for processes in
$\Dbb$ but then may have several different solutions as noted by
Majewski \cite{maj-fluid}. We give in the appendix a simple example.

To circumvent this difficulty, we adopt the following strategy. We
find a domain $\Dcal_\JN\subset\Ecal%\Dcal_\JN\subset \Dbb_0(\real^K_+)\times\Dbb_0(\Mbb^K)\times \Dbb(\real^K_+)
$ satisfying the following constraints:
\begin{itemize}
\item the map $\Psi$ is well defined on $\Dcal_\JN$;
\item any solution $(\A,\D)$ of the fixed point equation
(\ref{fixpointeq}) associated with a "continuous" Jackson network
$\JN$ can be approximated by a sequence $\{\JN_n\}\in
\Dcal_\JN^\nat$ such that
\begin{eqnarray}
\label{a}\JN_n&\rightarrow& \JN,\\
\label{b}\Psi(\JN_n)&\rightarrow& (\A,\D),\\
\label{c}\Isf^\JN(\JN_n)&\rightarrow& \Isf^\JN(\JN).
\end{eqnarray}
\end{itemize}

Hence in order to remove the quote from (\ref{eq:quote}), we follow
a quite standard method of proofs for large deviations of stochastic
processes analogue with the theory of weak convergence
\cite{fenkur}: it consists of first verifying a compactness
condition and then showing that there is only one possible limit. In
our context, we proceed as follows:
\begin{enumerate}
\item we show that our sequence of processes is exponentially tight;
\item we use $\Dcal_\JN$ to determine the rate function.
\end{enumerate}

In Section \ref{sec:contr}, we give the theoretical framework that
shows how any domain verifying assumptions (\ref{a}), (\ref{b}) and
(\ref{c}) determines the rate function. This result is stated in
great generality (without any reference to our specific problem) and
could be of independent interest since this method of proof could be
applied to other dynamical systems (with discontinuous statistics).

\subsection*{Notation}
For $(E,d,\leq)$ a complete, separable metric space with partial
order $\leq$, we denote by $\Dbb(E)$ the space of cadlag
non-decreasing $E$-valued functions defined on $\real_+$ with
Skorohod ($J_1$) topology and by $\Cbb(E)$ the space of continuous
non-decreasing $E$-valued functions defined on $\real_+$. Restricted
to $\Cbb(E)$ the Skorohod topology is just the compact uniform
topology.
%For any non-decreasing function $f$, we denote $f^{\leftarrow}(x)=\inf\{t,\:f(t)\geq x\}$ the pseudo-inverse of $f$.

For $x,y\in \real^K$, we write $x\leq y$ if $x^{(i)}\leq y^{(i)}$
for all $i$. We denote by $\wedge$ the minimum and by $\vee$ the
maximum in $\real^K$. For $\Xbf, \Ybf \in \Dbb(\real_+^K)$, we write
$\Xbf\leq \Ybf$ if $\Xbf(t)\leq \Ybf(t)$ for all $t\geq 0$ and for
maps $F,G \in \Dbb(\real_+^K)^2$, we denote $F\leq G$ if
$F(\Xbf)\leq G(\Xbf)$ for all $\Xbf\in \Dbb(\real_+^K)$. For $x\in
\real^K$, we denote $\|x\|=\vee_{i=1}^K x^{(i)}$ and for
$\Xbf\in\Dbb(\real_+^K)$, we denote $\|\Xbf\| = \sup_t \|\Xbf(t)\|$.
We denote $\Dbb_0(E)=\{f\in \Dbb(E),\:f(0)=0\}$ and
$\Cbb_0(E)=\{f\in \Cbb(E),\:f(0)=0\}$.

A piecewise linear function is a continuous function such that there
exists a partition $\tau=(t_0=0<t_1<\dots)$ with $t_k\rightarrow
\infty$ and such that the function is linear on each interval
$(t_k,t_{k+1})$. For any function $f\in \Dbb(\real_+^K)$, we define
the polygonal approximation of $f$ with step $1/n$ as the (piecewise
linear) function
\begin{eqnarray*}
f_n(t) = f\left(\frac{\lfloor nt\rfloor}{n}\right)+\left(nt- \lfloor
nt\rfloor\right)\left(f\left(\frac{\lfloor
nt\rfloor+1}{n}\right)-f\left(\frac{\lfloor
nt\rfloor}{n}\right)\right)
\end{eqnarray*}

$\Mbb^K$ is the set of substochastic matrices of size $K\times K$.
For $M\in \Mbb^K$, we denote by $\rho(M)$ its spectral radius, by
$M^t$ its transpose and $M^{(i)}$ denotes the line $M^{(i)}=
(M^{(i,1)},\dots M^{(i,K)})$. In particular, we will identify a
function $\Pbf\in \Dbb(\Mbb^K)$ with its $K$ components
$\Pbf^{(i)}\in \Dbb(\real_+^K)$, where $\Pbf^{(i)}(t)=
(\Pbf^{(i,1)}(t),\dots \Pbf^{(i,K)}(t))$ with
$\sum_{j}\Pbf^{(i,j)}(t)\leq 1$ for all $t\geq 0$ and all $i$. Note
that for $M,N\in \Mbb^K$, we have $M\leq N$ if $M^{(i,j)}\leq
N^{(i,j)}$ for all $i$ and $j$.

We will use the Kullback-Leibler information divergence, which is a
nonsymmetric measure of distance between distributions in the sense
that for any two distributions $P$ and $R$ on $\Xcal^k$ where
$\Xcal$ is a finite set,
\begin{eqnarray*}
\Dsf(P\|R) = \sum_{x\in \Xcal^k}P(x)\log\left(
\frac{P(x)}{R(x)}\right),
\end{eqnarray*}
is nonnegative and equals $0$ if and only if $P=R$. We use the
standard notational conventions $\log 0=-\infty$,
$\log\frac{1}{0}=\infty$ and $0 \log0=0\log\frac{0}{0}=0$. For any
fixed $R$, the divergence $\Dsf(P\|R)$ is a continuous function of
$P$ restricted to $\{P,\:S(P)\subset S(Q)\}$ where $S(P)$ denotes
the support of $P$ (see \cite{ccc}).

For $P\in \Mbb^K$, we denote by $\tilde{P}$ the $K\times (K+1)$
stochastic matrix obtained as follows: for all $i,j\leq K$,
$\tilde{P}^{(i,j)}=P^{(i,j)}$ and
$\tilde{P}^{(i,K+1)}=1-\sum_{k=1}^K P^{(i,k)}$. For $P,R\in \Mbb^K$,
we will denote
\begin{eqnarray*}
\tilde{\Dsf}(P\|R)&:=& \Dsf(\tilde{P}\|\tilde{R})\\
&=& \sum_{i,j\leq K} P^{(i,j)}\log\left(
\frac{P^{(i,j)}}{R^{i,j)}}\right)+\sum_{i\leq K}\left(1-\sum_k
P^{(i,k)}\right)\log\left(\frac{1-\sum_k P^{(i,k)}}{1-\sum_k
R^{(i,k)}} \right)\\
&=:& \sum _{i=1}^K \tilde{\Dsf}(P^{(i)}\|R^{(i)}).
\end{eqnarray*}

\section{An extension of the contraction principle}\label{sec:contr}

Let $\Ecal,\Fcal$ be complete separable metric spaces. Let
$G:\Ecal\times \Fcal\rightarrow \real$ be a continuous function. We
assume that there exists $\Dcal\subset \Ecal$, such that for all
$x\in \Dcal$, there exists an unique $y\in \Fcal$ such that
$G(x,y)=0$. We denote it by, $y=H(x)$ where $H:\Dcal\rightarrow
\Fcal$,
\begin{eqnarray*}
\forall x\in \Dcal,\quad G(x,y)=0 \:\Leftrightarrow \: y=H(x).
\end{eqnarray*}

\begin{prop}\label{prop:contraction}
Let $\{X_n\}_n$ be a sequence of $\Ecal$-valued random variables and
$\{Y_n\}_n$ be a sequence of $\Fcal$-valued random variables. We
assume that each sequence is exponentially tight. Assume that the
sequence $\{X_n\}_n$ satisfies a LDP with good rate function $I^X$
and that $G(X_n,Y_n)=0$  a.s. for all $n$.

We assume that for all $(x,y)$ such that $G(x,y)=0$ and
$\Isf^X(x)<\infty$, there exists a sequence \iffalse
$\{\Scal_n(x,y)\}_n\in \Dcal^\nat$ such that
\begin{itemize}
\item $\Scal_n(x,y) \rightarrow x$;
\item $H(\Scal_n(x,y)) \rightarrow y$;
\item $\Isf^X(\Scal_n(x,y))\rightarrow \Isf^X(x)$.
\end{itemize}
\fi $x_n\rightarrow x$, such that $x_n\in \Dcal$ for all $n$,
$H(x_n)\rightarrow y$ and $\Isf^X(x_n)\rightarrow \Isf^X(x)$. %For $(x,y)$ such that $G(x,y)=0$,
We denote by $\Scal(x,y)=\{x_n\}_n$ this sequence. If $G(x,y)\neq 0$
or $\Isf^X(x)=\infty$, we take $\Scal(x,y)=\emptyset$ and we denote
$\Scal(y) = \cup_x\{\Scal(x,y)\}$ (which might be empty).

Then the sequence $\{X_n,Y_n\}_n$ satisfies a LDP with good rate
function:
\begin{eqnarray}
\label{ratefunction}\Isf^{X,Y}(x,y):= \left\{\begin{array}{ll} \Isf^X(x), & G(x,y)=0,\\
\infty,& \mbox{otherwise.}
\end{array}\right.
\end{eqnarray}

In particular, if $X_n\in \Dcal$ for all $n$ and if the sequence
$\{H(X_n)\}_n$ is exponentially tight, then it satisfies a LDP in
$\Fcal$ with good rate function:
\begin{eqnarray}
\label{eq:rateH}\Isf^{H(X)}(y):=  \inf \{ \lim_{n\rightarrow
\infty}\Isf^X(x_n),\: \{x_n\}_n\in \Scal(y)\}.
\end{eqnarray}
\end{prop}

\begin{remark}\label{rem1}\begin{itemize}
\item There are alternative ways of expressing the rate function,
\begin{eqnarray*}
\Isf^{H(X)}(y)= \inf\{ \Isf^X(x), \: y\in H^x\},
\end{eqnarray*}
where $H^x:=\{y\in \Fcal,\: \exists x_n\rightarrow x,\:
H(x_n)\rightarrow y\}$. $\Isf^{H(X)}$ is the lower semicontinuous
regularization of the following function defined for $y\in
H(\Dcal)\subset \Fcal$,
\begin{eqnarray*}
\tilde{\Isf}^{H(X)}(y):=\inf\{ \Isf^X(x), \: y=H(x)\}.
\end{eqnarray*}
The main interest of the definition (\ref{eq:rateH}) is that the
rate function is computed only thanks to the sequences
$\Scal(x,y)\in \Dcal^\nat$.
\item Note that if $H(\Dcal)$ is closed (in particular if
$\Dcal=\Ecal$) then this proposition follows from the contraction
principle (for an extensive discussion of this principle, see the
work of Garcia \cite{garcia}). Roughly speaking, Proposition
\ref{prop:contraction} tells us that if $\Dcal$ is dense in a
certain sense in $\Ecal$, then the contraction principle still holds
for the map $H$.
\end{itemize}
\end{remark}
\proof Thanks to Lemma 3.6 of \cite{fenkur}, the sequence
$\{X_n,Y_n\}_n$ is exponentially tight. Then by Theorem 3.7 of
\cite{fenkur}, there exists a subsequence $\{n_k\}$ along which the
sequence $\{X_{n_k},Y_{n_k}\}_{n_k}$ satisfies a LDP with a good
rate function. If we can prove that there is a unique possible rate
function (that does not depend on the subsequence $\{n_k\}$) then
the proposition will follow.

Hence, for simplicity of notations, we still denote the extracted
subsequence by $\{X_n,Y_n\}_n$ and we assume that $\{X_n,Y_n\}_n$
satisfies a LDP with good rate function $\tilde{\Isf}^{X,Y}$. We
will show that $\tilde{\Isf}^{X,Y}=\Isf^{X,Y}$ given by
(\ref{ratefunction}).

Consider the continuous mappings $H_1$ and $H_2$ from $\Ecal\times
\Fcal$ to $\Ecal\times \Fcal \times \real$,
\begin{eqnarray*}
H_1(x,y):=(x,y,G(x,y)), \quad H_2(x,y):=(x,y,0).
\end{eqnarray*}

We have clearly $H_1(X_n,Y_n)=H_2(X_n,Y_n)$ a.s. Moreover thanks to
the contraction principle, $\{H_1(X_n,Y_n)\}_n$ and
$\{H_2(X_n,Y_n)\}_n$ satisfy LDPs with the good rate functions
\begin{eqnarray*}
\Isf^{H_1}(x,y,z) =\inf\{\tilde{\Isf}^{X,Y}(x,y), z=G(x,y)\}\quad
\Isf^{H_2}(x,y,z) =\inf\{\tilde{\Isf}^{X,Y}(x,y), z=0\},
\end{eqnarray*}
where $\inf\emptyset=\infty$. Since $H_1(X_n,Y_n)=H_2(X_n,Y_n)$, we
have $\Isf^{H_1}=\Isf^{H_2}$. Now we have,
\begin{eqnarray*}
\tilde{\Isf}^{X,Y}(x,y)=\inf_z\{\Isf^{H_1}(x,y,z)\}=\inf\{\tilde{\Isf}^{X,Y}(x,y),
G(x,y)=0\},
\end{eqnarray*}
hence $\tilde{\Isf}^{X,Y}(x,y)=\infty$ as soon as $G(x,y)\neq 0$. It
remains to show that $G(x,y)=0$ implies $\tilde{\Isf}^{X,Y}(x,y) =
\Isf^X(x)$. We have clearly $\Isf^X(x)\leq \tilde{\Isf}^{X,Y}(x,y)$
for all $(x,y)$ since $\{X_n\}$ satisfies a LDP with good rate
function
\begin{eqnarray*}
\Isf^{X}(x)=\inf\{\tilde{\Isf}^{X,Y}(x,y),\:y\in \Fcal,
\:G(x,y)=0\}.
\end{eqnarray*}
%In particular $\{Y_n\}$ satisfies a LDP with good rate function \begin{eqnarray*}I^{Y}(y):=\inf\{I^{X,Y}(x,y),\:x\in E, \:G(x,y)=0\}\end{eqnarray*}

In particular, the definition of $\Dcal$ implies
$\Isf^X(x)=\tilde{\Isf}^{X,Y}(x,H(x))$ for $x\in \Dcal$.

\iffalse For all $y\in F$, we assume that if $x$ is such that
$G(x,y)=0$ and $I^X(x)<\infty$, then there exists $x_n\rightarrow
x$, such that $x_n\in D$ for all $n$, $H(x_n)\rightarrow y$ and
$I^X(x_n)\rightarrow I^X(x)$. We denote
\begin{eqnarray*}
\tilde{I}^{X,Y}(x,y):= \left\{\begin{array}{ll} I^X(x), & G(x,y)=0,\\
\infty,& \mbox{otherwise.}
\end{array}\right.
\end{eqnarray*}

Since $I^{X}(x)\leq \tilde{I}^{X,Y}(x,y)$, we have
$I^{X,Y}(x,y)\leq\tilde{I}^{X,Y}(x,y)$.\fi

Take $(x,y)$ such that $G(x,y)=0$ and $\Isf^{X}(x)<\infty$. There
exists $x^*_n\rightarrow x$ with $x^*_n\in \Dcal$,
$H(x^*_n)\rightarrow y$ and $\Isf^X(x^*_n)\rightarrow \Isf^X(x)$.
Thanks to the lower semicontinuity property of $\tilde{\Isf}^{X,Y}$,
we can find for any $\delta>0$, an $\epsilon>0$ such that
\begin{eqnarray*}
\frac{1}{\delta}\wedge
\left(\tilde{\Isf}^{X,Y}(x,y)-\delta\right)\leq \inf_{z\in
B(y,\epsilon)}\tilde{\Isf}^{X,Y}(x,z),
\end{eqnarray*}
where $B(y,\epsilon)$ is the closed ball in $\Fcal$ of center $y$
and radius $\epsilon$.

Thanks to the lower semicontinuity of the function
$x\mapsto\inf_{z\in B(y,\epsilon)}\tilde{\Isf}^{X,Y}(x,z)$, we have
\begin{eqnarray*}
\inf_{z\in B(y,\epsilon)}\tilde{\Isf}^{X,Y}(x,z)&\leq&
\liminf_{x_n\rightarrow
x}\inf_{z\in B(y,\epsilon)}\tilde{\Isf}^{X,Y}(x_n,z)\\
&\leq& \liminf_{n\rightarrow \infty}\inf_{z\in
B(y,\epsilon)}\tilde{\Isf}^{X,Y}(x^*_n,z)\\
&\leq&\lim_{n\rightarrow \infty}\Isf^X(x^*_n)=\Isf^X(x),
\end{eqnarray*}
because $H(x^*_n)\in B(y,\epsilon)$ for sufficiently large $n$.
Hence we proved that for any $\delta>0$, $\frac{1}{\delta}\wedge
\left(\tilde{\Isf}^{X,Y}(x,y)-\delta\right)\leq \Isf^X(x)$ for
$(x,y)$ such that $G(x,y)=0$ and $\Isf^{X}(x)<\infty$, this
concludes the proof of (\ref{ratefunction}).

The various expressions of $\Isf^{H(X)}$ are now quite easy to
obtain from
\begin{eqnarray}
\label{eq:rateH1}\Isf^{H(X)}(y) = \inf\{ \Isf^X(x), \: G(x,y)=0\}.
\end{eqnarray}
For (\ref{eq:rateH}), note that since the set $\{x,\: G(x,y) = 0\}$
is closed the minimum in (\ref{eq:rateH1}) (if it is finite) is
attained for a certain $x^*$ with $G(x^*,y)=0$ and
$\Isf^X(x^*)<\infty$.

We prove now that
\begin{eqnarray*}
\inf\{\Isf^X(x),\:y\in H^x\}= \inf\{\Isf^X(x),\: G(x,y)=0\}.
\end{eqnarray*}
If $y\in H^x$, then there exists $x_n\rightarrow x$ such that
$H(x_n)\rightarrow y$. Hence by continuity of $G$, we have
$G(x,y)=0$. Now if $G(x,y)=0$ and $\Isf^X(x)<\infty$, it follows
from the assumptions that $y\in H^x$.

To see that the last expression in Remark \ref{rem1} is true, we
show that for any open set $O\subset \Fcal$, we have,
\begin{eqnarray}
\label{eq:ineqI}\inf_{y\in O}\Isf^{H(X)}(y) = \inf_{y\in O}
\{\Isf^X(x), \:y=H(x)\}.
\end{eqnarray}
For $y\in O$ and  any $x$ such that $G(x,y)=0$, there exists
$x_n\rightarrow x$, such that $H(x_n)\rightarrow y$ and
$\Isf^X(x_n)\rightarrow \Isf^X(x)$. Hence for $n$ sufficiently
large, we have $H(x_n)\in O$ and then
\begin{eqnarray*}
\inf_{y\in O} \{\Isf^X(x), \:y=H(x)\}\leq \inf_n \Isf^X(x_n)\leq
\Isf^X(x).
\end{eqnarray*}
Taking the minimum over all $x$ such that $G(x,y)=0$ gives the
$\geq$ inequality in (\ref{eq:ineqI}), the converse inequality is
obvious.
\endproof

\section{Queueing networks with Bernoulli routing: description and large deviations results}\label{sec:network}

\subsection{General setting and notation}

We start with the basic model for an isolated queue and refer to
\cite{ABR} for more details on the relationship with other models of
the literature.

The model for an isolated queue is in term of two primitive
quantities belonging to $\Dbb(\real_+)$: the arrival process $\A$
and the service process $\Sbf$. The departure process $\D$ is a
derived quantity that is obtained as a functional of the arrival and
service processes as follows:
\begin{eqnarray}
\label{ssq}\D(t) := \inf_{0\leq s\leq t}\left\{
\Sbf(t)-\Sbf(s)+\A(s)\right\}\wedge \Sbf(t).
\end{eqnarray}

From a mathematical point of view, if $\Rcal:\Dbb\rightarrow \Dbb$
(where $\Dbb$ is the space of cadlag $\real$-valued functions
defined on $\real_+$) is the one-dimensional Skorohod's reflection
map defined by $\Rcal(\Xbf)(t):=\sup_{0\leq s\leq
t}\left\{\Xbf(t)-\Xbf(s)\right\}\vee \Xbf(t)$. We have
$\D=\A-\Rcal(\A-\Sbf)$. It is easy to see that $\D\in \Dbb(\real_+)$
and $\D\leq \A$.

The queue length process is defined as the difference of the arrival
process and the departure process,
\begin{eqnarray*}
\Q(t):= \A(t)-\D(t)=\sup_{0\leq s\leq t}\left\{
\A(t)-\A(s)-(\Sbf(t)-\Sbf(s))\right\}\vee \left(
\A(t)-\Sbf(t)\right).
\end{eqnarray*}

If the arrival process $\A$ and the service process $\Sbf$ are
counting processes, this model is called a single queue with
autonomous service: the queue length is increased by one whenever
there is an arrival from the arrival process and the queue length is
decreased by one whenever there is an arrival from the service
process and the queue is not empty (see \cite{Borovkov}). Note in
particular that in the case where the process $\Sbf$ is a Poisson
point process, then this model is a standard $./M/1$ queue.

We now consider networks obtained by interconnecting queues modeled
by (\ref{ssq}) when the departure process of one queue is randomly
routed to the other queues as for Jackson networks. The networks we
consider are characterized by the fact that service times and
routing decisions are associated with stations and not with
customers. This means that we associate to each of the $K$ stations
three predefined counting processes: an arrival process, a service
process and a routing process. The arrival process and the service
process of station $k$ are described by the sequences of exogenous
arrival times $\{T^{(k)}_j\}_{j\geq 1}$ and service times
$\{\sigma^{(k)}_j\}_{j\geq 1}$. If there is no exogenous arrival at
station $k$, we use the convention $T^{(k)}_j=\infty$ for all $j$.
When the $j$-th customer has completed his service at station $k$,
he is sent to station $\nu^{(k)}_j$ (or leaves the network if
$\nu^{(k)}_j=K+1$) and is put at the end of the queue on this
station, where $\{\nu^{(k)}_j\}_{j\geq 1}$ is also a predefined
sequence, called the routing sequence. The sequences
$\{T^{(k)}_j\}_{j\geq 1}$, $\{\sigma^{(k)}_j\}_{j\geq 1}$ and
$\{\nu^{(k)}_j\}_{j\geq 1}$, where $k$ ranges over the set of
stations, are called the driving sequences of the network. A network
will be defined by $\left\{\{\sigma^{(k)}_j\}_{j\geq 1},\:
\{\nu^{(k)}_j\}_{j\geq 1},\:\{T^{(k)}_j\}_{j\geq 1},\: n^{(k)},\:
1\leq k\leq K\right\}$, where $(n^{(1)},\ldots,n^{(K)})$ describes
the initial condition. The interpretation is as follows: at time
$t=0$, in node $k$, there are $n^{(k)}$ customers with service times
$\sigma^{(k)}_1,\dots, \sigma^{(k)}_{n^{(k)}}$ (if appropriate,
$\sigma^{(k)}_1$ may be interpreted as a residual service time). In
particular at time $0$, the total number of customers in the network
is $n^{(1,K)}=n^{(1)}+\dots n^{(K)}$.

In what follows, we will describe the driving sequences thanks to
their associated counting functions. We will use the following
notation: $\sigma^{(k)}(1,n) = \sum_{j=1}^n\sigma^{(k)}_j$, for
$0\leq k\leq K$.

We define the sequence of networks ${\JN}_n=\{\Sbf_n(t),
\Pbf_n(t),\Nbf_n(t)\}$ with
\begin{eqnarray*}
\Nbf^{(i)}_n(t) &=&\frac{1}{n}\left(n^{(i)}_n+\sum_k \ind_{\{T^{(i)}_k\leq nt\}}\right),\\
\Sbf^{(i)}_n(t) &=& \frac{1}{n} \sum_k
\ind_{\{\sigma^{(i)}(1,k)\leq nt\}},\\
\Pbf_n^{(i,j)}(t) &=& \frac{1}{n} \sum_{k\leq
nt}\ind_{\{\nu^{(i)}_k=j\}}.
\end{eqnarray*}

Note that we allow the initial queue length to depend on $n$,
$\Nbf^{(i)}_n(0)=n^{(i)}_n$ but the other driving sequences
describing the arrival times, the service times and the routing
decisions do not depend on $n$. Note also that if there is no
exogenous arrival at station $i$, we have
$\Nbf^{(i)}_n(t)=\Nbf^{(i)}_n(0)$ for all $t\geq 0$.

For the network $\JN_n$, we denote the corresponding input and
output processes of each queue $k$ of the network by $\A_n^{(k)}$
and $\D_n^{(k)}$ respectively. We will use the following notation
$\A_n =(\A_n^{(1)},\dots, \A_n^{(K)})$ and $\D_n=(\D_n^{(1)},\dots,
\D_n^{(K)})$. We now describe how the processes $\A_n$ and $\D_n$
are obtained form $\JN_n$.

We define the map $\Gamma:\Dbb_0(\real^K_+)\times \Dbb_0(
\Mbb^K)\times \Dbb(\real_+^K)\rightarrow\Dbb(\real^K_+)$ as follows:
\begin{eqnarray*}
\Gamma(\Xbf,\Pbf,\Nbf)^{(i)}(t) := \Nbf^{(i)}(t) +
\sum_{j=1}^K\Pbf^{(j,i)}(\Xbf^{(j)}(t)).
\end{eqnarray*}

The following lemma is straightforward.
\begin{lem}
\label{lem:contGamma}The map $\Gamma$ is continuous for the compact uniform topology and non-decreasing in its first argument.% on $\Cbb(\real^K_+)\times \Cbb( \Mbb^K)\times \real_+^K$.
\end{lem}

We define the map
$\Phi:\Dbb(\real^K_+)\times\Dbb_0(\real^K_+)\rightarrow\Dbb_0(\real^K_+)$
as follows:
\begin{eqnarray*}
\Phi(\Xbf,\Ybf)^{(i)}(t):=  \inf_{0\leq s\leq t}\left\{
\Ybf^{(i)}(t)-\Ybf^{(i)}(s)+\Xbf^{(i)}(s)\right\}\wedge
\Ybf^{(i)}(t).
\end{eqnarray*}

\begin{lem}
\label{lem:contPhi} The map $\Phi$ is continuous for the compact uniform topology and non-decreasing in its first argument.% on $\Dbb(\real^K_+)\times\Dbb(\real^K_+)$.
\end{lem}
\proof We can clearly consider the map $\Phi$ with $K=1$ only. Let
$\Rcal$ be the one-dimensional reflection map, we have
$\Phi(\Xbf,\Ybf)=\Xbf-\Rcal(\Xbf-\Ybf)$. It is easy to see that for
any $T>0$,
\begin{eqnarray*}
\sup_{0\leq t\leq T}|\Rcal(\Xbf)(t) -\Rcal(\Xbf')(t)|\leq 2
\sup_{0\leq t\leq T}|\Xbf(t) -\Xbf'(t)|,
\end{eqnarray*}
from which the continuity of $\Phi$ follows. Its monotonicity is
obvious.
\endproof
\begin{remark}\label{rem:conc}Consider the mapping $\Phi$ with
$K=1$ and $\Ybf(t)=\mu t$, with $\mu\geq 0$. If $\mu=0$, since
$\Phi(\Xbf,\Ybf)\leq \Ybf$, we have $\Phi(\Xbf,\Ybf)(t)=0$ for all
$t$. If $\mu\not= 0$, we have $\Phi(\Xbf,\Ybf)(t) = \inf_{0\leq
s\leq t} \{\Xbf(s)+\mu(t-s)\}$. Moreover if $\Xbf$ is a concave
function, then this equation reduces to $\Phi(\Xbf,\Ybf)(t)=
\Xbf(t)\wedge \mu t$. Hence we can write
\begin{eqnarray*}
\mbox{$\Ybf(t)=\mu t$, with $\mu\geq 0$} &\Rightarrow &\Phi(\Xbf,\Ybf)(t) = \mu t \wedge \inf_{0\leq s\leq t} \{\Xbf(s)+\mu(t-s)\},\\
\mbox{if moreover $\Xbf$ is a concave function } & \Rightarrow &
\Phi(\Xbf,\Ybf)(t) = \mu t \wedge \Xbf(t).
\end{eqnarray*}
\end{remark}

It is easy to adapt the proof of Theorem 2.1 of \cite{cm:disc} or
Proposition 2.1 of \cite{lel} to show that the following fixed-point
equation:
\begin{eqnarray}
\label{fixpointeq1}
\left\{\begin{array}{lcl}\A_n &=& \Gamma(\D_n,\Pbf_n,\Nbf_n )=\Gamma(\D_n,\JN_n),\\
\D_n &=& \Phi(\A_n,\Sbf_n)=\Phi(\A_n,\JN_n),\end{array}\right.
\end{eqnarray}
has an unique solution when each component of $n\Sbf_n$, $n\Pbf_n$
and $n\Nbf_n$ is a counting function (i.e. non-decreasing function
of $\Dbb(\real^K_+)$ or $\Dbb(\Mbb^K)$ that is piece-wise constant
with jumps of size one). In this case the corresponding functions
$n\A_n$ and $n\D_n$ are also counting functions and we denote the
solution of (\ref{fixpointeq1}) by
$\Psi(\Sbf_n,\Pbf_n,\Nbf_n)=\Psi(\JN_n)$.

\begin{remark}
Note that the only difference between our model and generalized
Jackson networks as described in \cite{lel} resides in the queueing
mechanism (\ref{ssq}) which is sometimes called autonomous. Consider
a network $\JN=\{\Sbf,\Pbf, \Nbf\}$ where the processes are counting
processes. Then due to some monotonicity arguments, it is possible
to
relate (see \cite{CTK}):\\
- the process $(\tilde{\A},\tilde{\D})$ associated to $\JN$ with
the dynamic described in \cite{lel};\\
- the processes $\Psi(\JN)=(\A,\D)$ solution of the fixed point
equation.\\
Note that in the case where the process $\Sbf$ is a
Poisson point process, our model is exactly a Jackson network (see
\cite{Apq}).

\iffalse We will see in the next section that our stochastic
assumptions allow us to cover this case. Hence our large deviation
result is a generalization of \cite{atardupuis} and \cite{ign}. \fi
\end{remark}

\subsection{Stochastic assumptions}\label{sec:sto}

In what follows, it will be important to distinguish the nodes of
the network that do not receive any exogenous customer, i.e. the
nodes $i\in \Scal^c$ with $\Scal=\{i,\:T^{(i)}_1\leq \infty \}$. A
network $\JN=\{\Sbf,\Pbf,\Nbf\}$ is an object in $\Ecal\subset
\Dbb_0(\real_+^K) \times \Dbb_0(\Mbb^K)\times \Dbb(\real_+^K)$, with
the additional constraints: \begin{enumerate}
\item $\Nbf^{(i)}(t)=\Nbf^{(i)}(0)$ for all $t$, for $i\notin
\Scal$;
\item for all $0\leq v\leq u$, we have $\sum_{j=1}^K \Pbf^{(i,j)}(u)-\Pbf^{(i,j)}(v)\leq (u-v)$.%, in particular, we have $\dot{\Pbf}\in \Dbb_0(\Mbb^K)$ if it exists.
\end{enumerate}
Note that $\Ecal$ is closed in $\Dbb_0(\real_+^K) \times
\Dbb_0(\Mbb^K)\times \Dbb(\real_+^K)$.

We define for $(s^{(1)},\dots, s^{(K)})\in \real_+^K$ and
$(n^{(1)},\dots, n^{(K)})\in \real_+^K$, the functions
%$\Isf^\Nbf(n^{(1)},\dots, n^{(K)})$ and
\begin{eqnarray*}
\Isf^\Sbf(s^{(1)},\dots, s^{(K)}) &=& \sum_{i=1}^K
\Isf^{\Sbf^{(i)}}(s^{(i)}),\\
\Isf^\Nbf(n^{(1)},\dots, n^{(K)}) &=& \sum_{i\in \Scal}
\Isf^{\Nbf^{(i)}}(n^{(i)})+\infty\ind_{\{n^{(i)}>0, \:i\notin
\Scal\}},
\end{eqnarray*}
where each $\Isf^{\Sbf^{(i)}}$%:\real_+\mapsto \real_+\cup\{+\infty\}$
(resp. $\Isf^{\Nbf^{(i)}}$ for $i\in \Scal$) is a
$[0,\infty]$-valued convex good rate function, attaining zero on
$\real_+$ admitting a unique minimum at the point $\mu^{(i)}$ (resp.
$\lambda^{(i)}$ for $i\in \Scal$) and with a domain open on the
right.

We assume that the sequence ${\JN}_n=\{\Sbf_n(t),
\Pbf_n(t),\Nbf_n(t)\}$ satisfies a LDP in the space
%$\Dbb_0(\real_+^K)\times \Dbb_0(\Mbb^K)\times \Dbb(\real_+^K)$
$\Ecal$ with a good rate function $\Isf^{\JN}$ given by
\begin{eqnarray}
\label{ass:rate}\Isf^{\JN}(\Sbf, \Pbf,\Nbf) :=
\Isf^0(\Nbf(0))+\int_0^\infty
\Isf^\Sbf(\dot{\Sbf}(t))+\tilde{\Dsf}(\dot{\Pbf}(t)\|R)+\Isf^\Nbf(\dot{\Nbf}(t))
dt,
\end{eqnarray}
if the argument functions are absolutely continuous and equal to
infinity otherwise. %By definition $I^\JN$ is a good rate function.

We make the following assumptions on the matrix $R$:
\begin{enumerate}
\item We assume that $\rho(R)<1$.
\item We assume that for all $1\leq i\leq K$, we have
\begin{eqnarray}
\label{eq:RO}(\Ncal+\Ncal R+\dots+\Ncal R^K)^{(i)}>0,
\end{eqnarray}
where $\Ncal$ is the line vector of $\real_+^K$ defined by
$\Ncal^{(i)}=\ind_{\{i\in \Scal\}}$.
\end{enumerate}

We show now that our stochastic assumptions cover the case where
$\Sbf^{(k)}$ and $\Nbf^{(i)}$ (with $i\in \Scal$) are independent
and correspond to renewal processes and where the routing is a
Bernouilli routing associated with the matrix $R$ that satisfies
previous assumption.

We recall here some results of Puhalskii \cite{puh:ssq} concerning
large deviations of renewal processes and show that our assumptions
on the rate function (\ref{ass:rate}) are satisfied in the i.i.d
case. Denote by $\{\zeta_i,\:i\geq 1\}$ a sequence of non-negative
i.i.d. random variables with positive mean. Let
\begin{eqnarray}
\nonumber \alpha(\theta) &=& \log \esp\left[ e^{\theta \zeta_1}\right],\\
\nonumber \theta^* &=& \sup\{\theta>0, \:\alpha(\theta)<\infty \},\\
\label{eq:equal}\alpha^*(x) &=& \sup_{\theta}\{\theta x-\alpha(\theta)\}=\sup_{\theta< \theta^*}\{\theta x-\alpha(\theta)\},\\
\nonumber g(x)&=& x\alpha^*(1/x)=\sup_{\theta<\theta^*}\{\theta
-x\alpha(\theta)\}.
\end{eqnarray}
Note that the function $\alpha$ is a convex function and
differentiable on $(-\infty,\theta^*)$ with
$\alpha'(0)=\esp[\zeta_1]>0$. In particular, we have
$\lim_{\theta\uparrow\theta^*}\alpha(\theta)=\infty$, from which we
get the equality in (\ref{eq:equal}). The functions $\alpha^*$ and
$g$ are convex rate functions.
%Moreover if $\theta^*>0$ then $\alpha^*$ is a good rate function.
Introduce the sequence of
processes $\{\Cbf_n\}_n$:
\begin{eqnarray*}
\Cbf_n(t) = \frac{1}{n} \sum_{i}\ind_{\{\sum_{j=1}^i\zeta_j\leq
nt\}}.
\end{eqnarray*}
Then Theorem 3.1 of \cite{puh:ssq} gives: If $\prob(\zeta_1>0)=1$,
then the sequence $\{\Cbf_n\}_n$ satisfies a LDP in $\Dbb(\real_+)$
with the good rate function
\begin{eqnarray*}
\Isf^\Cbf(\x) =\left\{\begin{array}{ll} \int_0^\infty
g(\dot{\x}(t))dt,
&\mbox{if $\x\in \Cbb(\real_+)$ is absolutely continuous,}\\
\infty,&\mbox{otherwise.}
\end{array}\right.
\end{eqnarray*}
It then follows that $g$ is a good rate function. Moreover, we have
$\essinf \zeta_1 =0$ if and only if $g(x)$ is finite for all $x\geq
\esp[\zeta_1]^{-1}$ (note in particular, that in this case, the
domain of $g$ is open on the right). The proof of this fact can be
found in \cite{Puh:jack} and follows the argument: from
$\alpha(\theta) \geq \esp[\zeta_1]\theta$, we have
$g\left(\esp[\zeta_1]^{-1}\right)=0$ and for all $x\geq
\esp[\zeta_1]^{-1}$, we have $g(x) = \sup_{\theta\leq
0}\left\{\theta -x\alpha(\theta) \right\}$. If $\essinf \zeta_1=0$,
we have for arbitrary $\epsilon>0$ and for $\theta\leq 0$,
\begin{eqnarray*}
\alpha(\theta) = \log\esp\left[ e^{\theta \zeta_1}\right]\geq \theta
\epsilon+\log \prob(\zeta_1<\epsilon),
\end{eqnarray*}
hence for $x>\epsilon^{-1}$, we have $g(x)\leq -\log
\prob(\zeta_1<\epsilon)$. It is clear that if $\essinf \zeta_1>0$,
then for any $x>\essinf \zeta_1^{-1}$, we have $g(x)=\infty$.

Concerning the large deviations of the routing processes given in
term of the Kullback-Leibler information divergence, it follows
directly from Corollary 6.1 of \cite{funcldp} in the case of
Bernouilli routing, i.e. when the sequences $\{\nu^{(k)}_j\}_{j\geq
1}$ are sequences of i.i.d. random variables in $[1,K]$ and
independent in $k$ such that
\begin{eqnarray*}
\prob(\nu^{(k)}_1=i)=R^{(k,i)}.
\end{eqnarray*}

%\newpage

\subsection{Sample path large deviations for the queue length
process}\label{sec:result}

We now return to the sequence of queueing networks defined in
Section \ref{sec:network}. Recall that $(\A_n,\D_n)$ correspond to
the arrival and departure processes from each station. We now give
our theorem for the queue length process defined as
$\Q_n=\A_n-\D_n$.

\begin{theorem}\label{th:qu}
The sequence of processes $\{\Q_n\}_n$ satisfies a LDP in
$\Dbb(\real^K_+)$ with good rate function that is finite for $\Q$
absolutely continuous given by:
\begin{eqnarray*}
\Isf^{0}(\Q(0))+\Isf^\Q_{\Q(0)}(\Q),
\end{eqnarray*}
where for $q\geq 0$, $\Isf^\Q_q(.)$ is a good rate function that is
finite for absolutely continuous $\Q$ such that $\Q(0)=q$ and given
by:
\begin{eqnarray*}
\Isf^\Q_q(\Q) := \int_0^\infty \Hsf^\Q(\Q(s),\dot{\Q}(s))ds,
\end{eqnarray*}
where $\Hsf^\Q$ is given by,
\begin{eqnarray*}
\Hsf^\Q(Q,\dot{Q}):=\inf\left\{\sum_{i\in
E(Q)}\Isf^{\Sbf^{(i)}}(D^{(i)})\ind_{\{D^{(i)}>\mu^{(i)}\}}+\sum_{i\notin
E(Q)}\Isf^{\Sbf^{(i)}}(D^{(i)}) +\sum_i D^{(i)}
\tilde{\Dsf}(P^{(i)}\|R^{(i)})+\Isf^\Nbf(N)\right\}
\end{eqnarray*}
where $E(Q)=\{i,\: Q^{(i)}=0\}$ and the infimum is taken over the
set of $(D,P,N)\in \real_+^K\times \Mbb^K\times \real_+^K$ such that
\begin{eqnarray*}
\dot{Q} = N+(P^t-Id)D.%\quad \mbox{and,}\quad \rho(P)<1.
\end{eqnarray*}
\end{theorem}

In \cite{Puh:jack}, Puhalskii obtains a LDP for the queue length
process of a generalized Jackson network with a rate function that
coincides with Theorem \ref{th:qu}. Note that our model is slightly
different here since we model the dynamic of a queue by a reflection
mapping. Still in the case of Poisson processes for the inputs, both
models correspond to the (exponential distribution) Jackson network.
Recall that the rate function for a Poisson process of rate
$\lambda$ is given by (we keep the same notation as in
\ref{sec:sto}),
\begin{eqnarray}\label{eq:pois}
\Isf^{\Cbf}(\xbf)=\int_0^\infty \lambda
\dot{\xbf}(t)\log\frac{\dot{\xbf}(t)}{\lambda}-\dot{\xbf}(t)+\lambda\:
dt,
\end{eqnarray}
for absolutely continuous functions $\xbf\in \Cbb(\real_+)$. Hence
if we replace (\ref{eq:pois}) in the expression of $I^\Q_q$, we
obtain the rate function for the large deviations of a Jackson
network. In this specific case, the rate function has been obtained
in different forms by Atar and Dupuis \cite{atardupuis} and
Igniatiouk-Robert \cite{ign} and some bounds have been computed by
Majewski \cite{maj-ldbound}. Compare to these results, our
representation has the advantage of being quite intuitive, in the
sense that each term is easy to interpret. If we interpret $D,P,N$
as instantaneous departure, routing and exogenous arrival rates,
then $N+(P^t-Id)D$ is just the vector of rates at which the queue
lengths vary. Hence given a rate of change of $\Q$, the system
behaves in such a way to minimize the instantaneous "costs" of
departure, routing and exogenous arrival rates over all the rates
that yield the desired $\dot{\Q}$.

From a methodological point of view, the argument of \cite{Puh:jack}
is quite different from ours since the density condition (that we
could compare to our Proposition \ref{prop:contraction}) is verified
on the rate function $\Isf^{\Q}_q$ (see condition (D) in
\cite{Puh:jack}) whereas we are checking the density argument on the
rate function of the inputs.

\section{Extension of $\Psi$ to piece-wise linear networks}\label{sec:extpsi}

In this section we consider processes that are continuous, i.e. in
$\Cbb(E)$, hence topological concepts refer to the compact uniform
topology.

We first recall Proposition 3.2 of \cite{lel},
\begin{proposition}
\label{prop:fixpoint:lel} Given a $K\times K$ substochastic matrix
$P$ with $\rho(P)<1$ and vectors $(\alpha,y)\in \real_+^{2K}$, the
fixed point equation
\begin{eqnarray*}
x^{(i)}&=& \alpha^{(i)}+\sum_{j=1}^K P^{(j,i)}\left(x^{(j)}\wedge
y^{(j)}\right),
\end{eqnarray*}
has a unique solution $x(y,P,\alpha)$. Moreover, $(y,\alpha)\mapsto
x(y,P,\alpha)$ is a continuous non-decreasing function.
\end{proposition}

We first consider a linear network $\JN$ and show that the mapping
$\Psi$ (defined as the solution of the fixed-point Equation
(\ref{fixpointeq1})) is well defined for such a network. By linear,
we mean the following $\Nbf^{(i)}(t)=N^{(i)}+\lambda^{(i)}t$, with
$\lambda^{(i)}\geq 0$ and $N^{(i)}\in \real_+$, $\Sbf^{(i)}(t)
=\mu^{(i)}t$, with $\mu^{(i)}\geq 0$, and $\Pbf^{(i,j)}(t) =
P^{(i,j)}t$. We assume that $\rho(P)<1$.

\begin{lem}
\label{lem:linear}Under previous assumptions, the fixed point
equation (\ref{fixpointeq}) has an unique solution
$\Xbf_f[\mu,P,N,\lambda](t)=x(\mu t,P,N+\lambda t)$, where $\mu
=(\mu^{(i)})_i$, $N=(N^{(i)})_i$ and $\lambda=(\lambda^{(i)})_i$.
\end{lem}
\proof Since $\mu,P,N,\lambda$ are fixed here, we omit to explicitly
write the dependence in these variables. In this case, the fixed
point equation (\ref{fixpointeq}) reduces to (see Remark
\ref{rem:conc})
\begin{eqnarray}
\label{fixpointeq:lin}\left\{\begin{array}{lcl}\A^{(i)}(t) &=& N^{(i)}+\lambda^{(i)}t+\sum_{j=1}^K P^{(j,i)}\D^{(j)}(t),\\
\D^{(i)}(t) &=& \mu^{(i)}t \wedge \inf_{0\leq s\leq t}\{
\A^{(i)}(s)+\mu^{(i)}(t-s)\}.%\leq \A^{(i)}(t) \wedge\mu^{(i)}t.
\end{array}\right.
\end{eqnarray}

Thanks to Proposition \ref{prop:fixpoint:lel}, $\Xbf_f(t)=x(\mu
t,P,N+\lambda t)$ is the unique solution of the fixed point equation
\begin{eqnarray}
\label{fixpointeq:lin2}\left\{\begin{array}{lcl}\A^{(i)}(t) &=& N^{(i)}+\lambda^{(i)}t+\sum_{j=1}^K P^{(j,i)}\D^{(j)}(t),\\
\D^{(i)}(t) &=&  \A^{(i)}(t) \wedge \mu^{(i)}t.\end{array}\right.
\end{eqnarray}
We prove now that $\Xbf_f$ is the unique solution of the fixed point
equation (\ref{fixpointeq:lin}).

For simplicity, we denote the fixed point equation
(\ref{fixpointeq:lin}), resp. (\ref{fixpointeq:lin2}), by
$\A=F(\A)$, resp. by $\A=\tilde{F}(\A)$. Note that these functions
are non-decreasing, continuous and such that $F\leq \tilde{F}$.

From $\0\leq \Xbf_f$, we get $\0\leq F(\0)\leq \tilde{F}(\0)\leq
\tilde{F}(\Xbf_f)$. Hence $F^n(\0)\nearrow \Lbf\leq \Xbf_f$ and
$F(\Lbf)=\Lbf$. Moreover for any solution $\Ybf$ of the fixed point
equation (\ref{fixpointeq:lin}), we have $\Lbf\leq \Ybf \leq \Xbf_f$
because $\Ybf=F(\Ybf)\leq \tilde{F}(\Ybf)$ and
$\tilde{F}^n(\Ybf)\nearrow \Xbf_f$.

Since $\0$ is a concave function, we have $F(\0)=\tilde{F}(\0)$ and
hence it is still a concave function. Hence we have
$\tilde{F}^n(\0)=F^n(\0)$ since the image by $\tilde{F}$ of a
concave function is a concave function and $F=\tilde{F}$ on the
subspace of concave functions. Hence we have $\Lbf = \Xbf_f$ which
concludes the proof.
\endproof

In order, to extend $\Psi$ to piece-wise linear networks, we proceed
step by step on each interval where the driving functions
$\Sbf,\Pbf,\Nbf$ are linear. The following lemma allows to glue the
constructed solution on each adjacent interval. \iffalse In a
queueing context, this lemma says that the output of a single server
queue fed by the arrival process $\A$ and service time process
$\Sbf$ viewed from time $u$ is just the same as the output process
of a single server queue that we start at time $u$ with arrival
process $\tilde{\A}(t)=\A(t+u)-\A(u)+\A(u)-\D(u)$ (i.e. with the
same increment as the original process on this period of time plus
an additional bulk corresponding to the queue length at time $u$)
and with service time process $\tilde{\Sbf}(t)$.\fi
\begin{lem}\label{lem:trans}
Let $\A,\Sbf\in \Dbb(\real_+)\times \Dbb_0(\real_+)$ and
$\D=\Phi(\A,\Sbf)$. Define $\tilde{\A},\tilde{\Sbf}\in
\Dbb(\real_+)\times \Dbb_0(\real_+)$ as follows
\begin{eqnarray*}
\tilde{\A}(t) &:=& \A(t+u)-\D(u),\\
\tilde{\Sbf}(t) &:=& \Sbf(t+u)-\Sbf(u).
\end{eqnarray*}
Let $\tilde{\D}=\Phi(\tilde{\A},\tilde{\Sbf})$, then we have
\begin{eqnarray*}
\tilde{\D}(t) = \D(t+u)-\D(u).
\end{eqnarray*}
\end{lem}
\proof We show that for $\D=\Phi(\A,\Sbf)$, we have
\begin{eqnarray*}
\D(t+u)-\D(u) = \inf_{u\leq s\leq t+u}\left\{ \Sbf(t+u)
-\Sbf(s)+\A(s)-\D(u)\right\}\wedge \left\{\Sbf(t+u)-\Sbf(u)
\right\},
\end{eqnarray*}
from which the lemma follows.

We write
\begin{eqnarray*}
\D(t+u)-\D(u) &=& \inf_{0\leq s\leq u}\left\{ \Sbf(t+u)
-\Sbf(s)+\A(s)-\D(u)\right\}\\
&&\wedge \inf_{u\leq s\leq t+u}\left\{ \Sbf(t+u)
-\Sbf(s)+\A(s)-\D(u)\right\}\wedge \left\{\Sbf(t+u)-\D(u) \right\},
\end{eqnarray*}
Since $\D(u)\leq \Sbf(u)$, we have to prove that
\begin{eqnarray*}
\Sbf(t+u)-\Sbf(u) \geq \inf_{0\leq s\leq u}\left\{ \Sbf(t+u)
-\Sbf(s)+\A(s)-\D(u)\right\}\wedge \left\{\Sbf(t+u)-\D(u) \right\}.
\end{eqnarray*}
This will follow from,
\begin{eqnarray*}
\inf_{0\leq s\leq u}\left\{ \Sbf(t+u) -\Sbf(s)+\A(s)-\D(u)\right\}
&=&\Sbf(t+u)-\Sbf(u) +\inf_{0\leq s\leq u}\left\{ \Sbf(u)
-\Sbf(s)+\A(s)\right\}-\D(u)\\
&\leq&\Sbf(t+u)-\Sbf(u).
\end{eqnarray*}
\endproof

We consider now piece-wise linear networks: the functions
$u\mapsto\Nbf^{(i)}(u), u\mapsto\Sbf^{(i)}(u)$ and $u\mapsto
\Pbf^{(i,j)}(u)$ are continuous piece-wise linear functions such
that $\Nbf^{(i)}(0)\in\real_+$ and $\Sbf^{(i)}(0)=\Pbf^{(i,j)}(0)=0$
and $\rho(\dot{\Pbf}(t))<1$ for all $t\geq 0$.

\begin{prop}\label{prop:JNlinear}
For a piece-wise linear network, there exists an unique solution of
the fixed point equation (\ref{fixpointeq}). We still denote by
$\Psi$ the mapping  that to any piece-wise linear network $\JN$
associates the corresponding couple $(\A,\D)$.
\end{prop}
\proof The existence is a direct consequence of monotonicity
properties and continuity of the maps $\Gamma$ and $\Phi$. We define
the sequence of processes $\{\A[k],\D[k]\}_{k\geq 0}$ with the
recurrence equation:
\begin{eqnarray*}
\left\{ \begin{array}{l}
\A[k+1] = \Gamma(\D[k],\JN),\\
\D[k+1] = \Phi(\A[k+1],\JN),
\end{array}\right.
\end{eqnarray*}
and with initial condition $\D[0] = \0$. By the monotonicity
properties of $\Phi$ and $\Gamma$, we have
\begin{eqnarray*}
\0\leq\A[1] &\Rightarrow& \Phi(\0,\JN)=\0=\D[0] \leq
\Phi(\A[1],\JN)=\D[1]\\
&\Rightarrow& \Gamma(\D[0],\JN)=\A[1]\leq \Gamma(\D[1],\JN)=\A[2],
\end{eqnarray*}
and the sequence $\{\A[k],\D[k]\}_{k\geq 0}$ is increasing. Note
that $\D[k]\leq \Sbf$ and hence the following limits are well
defined
\begin{eqnarray*}
\lim_{k\rightarrow \infty} \A[k] = \A\quad \mbox{and,}\quad
\lim_{k\rightarrow \infty} \D[k] = \D.
\end{eqnarray*}
Since $\Gamma$ and $\Phi$ are continuous, $(\A,\D)$ is a solution of
the fixed point equation (\ref{fixpointeq}).

We now prove uniqueness. First recall that we call $\alpha$, a
partition of $\real_+$, any increasing sequence of points
$\alpha=\{a_n\}_n$ with $a_0=0$ and $a_n\rightarrow \infty$. For two
partitions $\alpha=\{a_n\}_n$ and $\beta=\{b_n\}_n$, we say that
$\gamma=\{g_n\}_n$ is the union of $\alpha$ and $\beta$ if $\gamma$
is a partition such that for all $n$ there exists $m$ such that
either $g_n=a_m$ or $g_n=b_m$.

Let $\tau=\{t_n\}_n$ be the union of the partitions associated with
each function $\Sbf, \Pbf, \Nbf$. %Let $\tau^{(i)}=\{t^{(i)}_n\}_n$ be the partition defined by\begin{eqnarray*}t^{(i)}_n = \inf\{t,\: \Sbf^{(i)}(t)\geq t_n\}.\end{eqnarray*}
We define for $x\in \real_+$, $d(x,\tau)=\min_n\{t_n-x,\: t_n> x\}>
0$.
%Finally let $\upsilon=\{u_n\}_n$ be the union of $\tau$ and the different $\tau^{(i)}$.

Assume that we are given two solutions of the fixed point equation
(\ref{fixpointeq}): $(\A_1,\D_1)$  and $(\A_2,\D_2)$. First note
that thanks to Lemmas \ref{lem:wPhi} and \ref{lem:wGamma}, any
solution of (\ref{fixpointeq}) is absolutely continuous. Let
$z=\inf\{t, \A_1(t)\neq \A_2(t)\}$, in particular, we have
$\A_1(t)=\A_2(t)$ and $\D_1(t)=\D_2(t)$ for all $t\leq z$.

Define $u=\min_i d(\D_\bullet^{(i)}(z),\tau)\wedge d(z,\tau)>0$,
where the notation $_\bullet$ can be replaced either by $_1$ or by
$_2$. We have that for $t\in [0,u]$,
\begin{eqnarray*}
\tilde{\Sbf}^{(i)}(t)&:=&\Sbf^{(i)}(z+t)-\Sbf^{(i)}(z)=t\mu^{(i)}, \\
\tilde{\Pbf}^{(i,j)}(t)&:=&\Pbf^{(i,j)}(\D_\bullet^{(i)}(z)+t)-\Pbf^{(i,j)}(\D_\bullet^{(i)}(z))=t P^{(i,j)}, \\
\tilde{\Nbf}^{(i)}&:=&\Nbf^{(i)}(z+t)-\Nbf^{(i)}(z)+\A_\bullet^{(i)}(z)-\D_\bullet^{(i)}(z)=t\lambda^{(i)}+\A_\bullet^{(i)}(z)-\D_\bullet^{(i)}(z),
\end{eqnarray*}
Let
$\tilde{\A}(t)=\Xbf_f[\mu,P,\A_\bullet(z)-\D_\bullet(z),\lambda](t)$
be the unique solution associated to the infinite horizon linear
network defined above. The associated departure process is
$\tilde{\D}(t) = \tilde{\A}(t)\wedge\mu t$. Let $v=\inf\{t,\:
\inf_i\tilde{\D}^{(i)}(t)=u\}$, in particular since
$\tilde{\D}^{(i)}(t)\leq \mu^{(i)}t$, we have $v>0$. In view of
Lemma \ref{lem:trans}, we have for $t\in (0,v)$,
\begin{eqnarray*}
\A_\bullet(t+z) = \tilde{\A}(t)+\D(z),\quad \D_\bullet(t+z)=
\tilde{\D}(t)+\D(z)
\end{eqnarray*}
this contradicts the fact that $z<\infty$ and concludes the proof.
\endproof

Let $\Ecal \subset \Dbb_0(\real_+^K) \times \Dbb_0(\Mbb^K)\times
\Dbb(\real_+^K)$ as defined at the beginning of Section
\ref{sec:sto} and $\Fcal=\Dbb(\real_+^K)\times \Dbb_0(\real_+^K)$.

For $\JN\in \Ecal$ and $(\A,\D)\in \Fcal$, we define the function
\begin{eqnarray*}
G(\JN,\A,\D)= \|(\A-\Gamma(\D,\JN),\D-\Phi(\A,\JN))\|.
\end{eqnarray*}
The function $G$ is continuous and such that
\begin{eqnarray*}
G(\JN,\A,\D)= 0 \Leftrightarrow\left\{\begin{array}{lcl}\A &=& \Gamma(\D,\JN),\\
\D &=& \Phi(\A,\JN).\end{array}\right.
\end{eqnarray*}

Let $\Dcal_\JN$ be the subspace of $\Ecal$ of piecewise linear
networks: namely $\JN=(\Sbf,\Pbf,\Nbf)\in \Dcal_\JN$ if the
functions $u\mapsto\Nbf^{(i)}(u), u\mapsto\Sbf^{(i)}(u)$ and
$u\mapsto \Pbf^{(i,j)}(u)$ are piecewise linear non-decreasing
functions such that $\rho(\dot{\Pbf}(t))<1$ for all $t\geq 0$ and
$\Nbf^{(i)}=\0$ for $i\notin \Scal$. We denote
$\dot{\JN}=(\dot{\Sbf},\dot{\Pbf},\dot{\Nbf})$.

We proved that
\begin{eqnarray*}
\forall \JN\in \Dcal_\JN,\quad G(\JN,\A,\D)=0 \Leftrightarrow
(\A,\D)=\Psi(\JN),
\end{eqnarray*}
where $\Psi$ has been explicitly defined above. We are exactly in
the framework of Section \ref{sec:contr}. In the next section we
construct the mapping $\Scal: \Ecal\times \Fcal\rightarrow
\Dcal_\JN^\nat$.

\section{Sample path large deviations}\label{sec:ldp}

In order to simplify the notations, we assume that $\Nbf_n(0)=0$ for
all $n$. This condition can be weakened to the standard condition:
\begin{eqnarray*}
\lim_{n\rightarrow
\infty}\frac{1}{n}\log\prob(\Nbf_n(0)>\epsilon)=0,
\end{eqnarray*}
for all $\epsilon>0$. In this case, we have $\Isf^0(x)=\infty$ for
all $x\neq0$ and $\Isf^0(0)=0$.

It is possible to deal with the case where the initial condition
satisfies a LDP as assumed in Theorem \ref{th:qu} by using a
standard conditioning argument (as done in \cite{puh:ssq} for
example).

\subsection{Construction of the approximating sequence}

This section is devoted to the proof of the following proposition:

\begin{proposition}\label{prop:JNapprox}
We consider $\JN=(\Sbf,\Pbf,\Nbf)\in \Ecal$ such that
$I^\JN(\JN)<\infty$ and such that there exists $(\A,\D)\in \Fcal$
that satisfies the fixed point equation (\ref{fixpointeq}) given by,
\begin{eqnarray*}
\left\{\begin{array}{lcl}\A &=&\Gamma(\D,\JN),\\
\D &=& \Phi(\A,\JN).\end{array}\right.
\end{eqnarray*}
There exists a sequence $\{\JN_n\}_n=\Scal(\JN,\A,\D)$ such that
\begin{eqnarray}
\label{eq:S1}\JN_n&\in& \Dcal_\JN\quad \mbox{for all $n$;}\\
\label{eq:S2}\JN_n &\rightarrow & \JN;\\%\quad \mbox{and,}\quad \Nbf_n(0)=\Nbf(0);\\
\label{eq:S3}\Psi(\JN_n)&\rightarrow& (\A,\D);\\
\label{eq:S4}\Isf^\JN(\JN_n)&\rightarrow& \Isf^\JN(\JN).
\end{eqnarray}
\end{proposition}

First note that since $\Isf^\JN(\JN)<\infty$, each process
$\Sbf,\Pbf,\Nbf$ is absolutely continuous and $\dot{\JN}$ is
well-defined. Moreover thanks to Lemma \ref{lem:propPhi}, the
processes $\A$ and $\D$ are absolutely continuous too.

The idea to construct the sequence $\{\JN_n\}_n$ is to consider the
piecewise approximation of the fixed point equation
(\ref{fixpointeq}). %This has to be done with some care and is the main technical part of the argument.
First consider the routing
equation $\A=\Gamma(\D,\JN)$ for times $t$ such that $nt\in \nat$,
\begin{eqnarray*}
\underbrace{\A^{(i)}(t+1/n)-\A^{(i)}(t)}_{\Delta^{(i)}_n(\A)(t)}
&=&\underbrace{\Nbf^{(i)}(t+1/n)-\Nbf^{(i)}(t)}_
{\Delta^{(i)}_n(\Nbf)(t)}+\sum_{j=1}^K
\dot{\tilde{\Pbf}}_n^{(j,i)}(\D^{(j)}(t+))\underbrace{(\D^{(j)}(t+1/n)-\D^{(j)}(t))}_{\Delta^{(j)}_n(\D)(t)},
\end{eqnarray*}
where we define the piece-wise linear process
$\tilde{\Pbf}_n^{(j,i)}(t)$ as follows, for $s\in
(\D^{(j)}(t),\D^{(j)}(t+1/n))$,
\begin{eqnarray*}
\dot{\tilde{\Pbf}}_n^{(j,i)}(s)
:=\frac{\Pbf^{(j,i)}(\D^{(j)}(t+1/n))-\Pbf^{(j,i)}(\D^{(j)}(t))}{\D^{(j)}(t+1/n)-\D^{(j)}(t)},
\end{eqnarray*}
if $\D^{(j)}(t+1/n)\neq \D^{(j)}(t)$, and we take
$\dot{\tilde{\Pbf}}_n^{(j,i)}(\D^{(j)}(t))=0$ otherwise. In other
words, we have
\begin{eqnarray*}
\tilde{\Pbf}_n^{(j,i)}(\D^{(j)}(t+1/n))-\tilde{\Pbf}_n^{(j,i)}(\D^{(j)}(t))&=&\dot{\tilde{\Pbf}}_n^{(j,i)}(\D^{(j)}(t+))(\D^{(j)}(t+1/n)-\D^{(j)}(t))\\
&=&\Pbf^{(j,i)}(\D^{(j)}(t+1/n))-\Pbf^{(j,i)}(\D^{(j)}(t))
\end{eqnarray*}

Note that $\{\dot{\tilde{\Pbf}}_n^{(j,i)}(t)\}_{i,j}\in \Mbb^K$
since we have by the definition of $\Ecal$,
\begin{eqnarray*}
\sum_i\Pbf^{(j,i)}(\D^{(j)}(t+1/n))-\Pbf^{(j,i)}(\D^{(j)}(t))\leq\D^{(j)}(t+1/n)-\D^{(j)}(t),
\end{eqnarray*}
but the matrix $(\dot{\tilde{\Pbf}}_n^{(j,i)}(\D^{(j)}(t+)))_{i,j}$
may not be of spectral radius less than $1$.

To circumvent this difficulty, we modify slightly the processes as
follows, (the variables $\eta,\epsilon_n,\delta$ will be made
precise latter)
\begin{eqnarray}
\label{eq:approx}\Delta_n^{(i)}(\A)+\frac{\eta^{(i)}}{n}&=&\Delta_n^{(i)}(\Nbf)+\frac{\delta^{(i)}}{n}\\
\nonumber&&+\sum_{j=1}^K
\left((1-\epsilon_n^{(j)})\dot{\tilde{\Pbf}}_n^{(j,i)}+\epsilon_n^{(j)}
R^{(j,i)}
\right)\left(\Delta_n^{(j)}(\D)+\frac{\eta^{(j)}}{n}\right),
\end{eqnarray}
where we omit to write the time $t$ and use the simplified notation
$\dot{\tilde{\Pbf}}_n^{(j,i)}=\dot{\tilde{\Pbf}}_n^{(j,i)}(\D^{(j)}(t+))$.

We have to find $\eta,\epsilon_n,\delta$ such that (\ref{eq:approx})
holds with $\eta^{(i)},\epsilon^{(i)}_n,\delta^{(i)}$ non-negative
and $\delta^{(i)}=0$ for $i\not\in \Scal$. These constraints are
satisfied by the following choice: first take $\delta$ such that
$\delta^{(i)}>0$ for all $i\in \Scal$ and $\delta^{(i)}=0$ for
$i\not\in \Scal$. Let $\eta(\delta)=\eta$ be the unique solution in
$\real_+^K$ of the following equation (recall that $\rho(R)<1$),
\begin{eqnarray*}
\eta^{(i)} = \delta^{(i)}+\sum_{j=1}^K \eta^{(j)}R^{(j,i)}.
\end{eqnarray*}
Note that $\eta^{(i)}>0$ for all $i$ thanks to (\ref{eq:RO}).
Finally let define $\epsilon_n(\delta)=\epsilon_n$ as follows
$\epsilon_n^{(i)} =
\frac{\eta^{(i)}}{n\Delta_n^{(i)}(\D)+\eta^{(i)}}\in (0,1]$ (note
that $\epsilon^{(i)}_n=1$ if and only if $\Delta_n^{(i)}(\D)=0$).

It is easy to see that (\ref{eq:approx}) holds since we have
\begin{eqnarray*}
(1-\epsilon_n^{(j)})\left(\Delta_n^{(j)}(\D)+\frac{\eta^{(j)}}{n}\right)
= \Delta_n^{(j)}(\D),&\mbox{or,}&
\epsilon_n^{(j)}\left(\Delta_n^{(j)}(\D)+\frac{\eta^{(j)}}{n}\right)
= \frac{\eta^{(j)}}{n},
\end{eqnarray*}
which imply respectively that
\begin{eqnarray*}
\Delta_n^{(i)}(\A)&=&\Delta_n^{(i)}(\Nbf)+\sum_{j=1}^K
(1-\epsilon_n^{(j)})\dot{\tilde{\Pbf}}_n^{(j,i)}\left(\Delta_n^{(j)}(\D)+\frac{\eta^{(j)}}{n}\right)\quad\mbox{and,}\\
\frac{\eta^{(i)}}{n}&=&\frac{\delta^{(i)}}{n}+\sum_{j=1}^K
\epsilon_n^{(j)} R^{(j,i)}
\left(\Delta_n^{(j)}(\D)+\frac{\eta^{(j)}}{n}\right),
\end{eqnarray*}
and summing these two equalities gives (\ref{eq:approx}).

%We take the following notations $U_\delta(t) =t\eta(\delta)$, and $V_\delta(t) =t\delta$.

For $\delta$ fixed, we define for $s\in
(\D^{(j)}(t)+t\eta(\delta),\D^{(j)}(t+1/n)+(t+1/n)\eta(\delta))$,
\begin{eqnarray*}
\dot{\Pbf}_{n,\delta}^{(j,i)}(s)
=(1-\epsilon_n^{(j)})\dot{\tilde{\Pbf}}_n^{(j,i)}(\D^{(j)}(t+))+\epsilon_n^{(j)}
R^{(j,i)},
\end{eqnarray*}
where $\epsilon_n(\delta)$ is defined as above. In view of Lemma
\ref{lem:specrad}, the matrix $\dot{\Pbf}_{n,\delta}^{(j,i)}(s)$ is
of spectral radius less than one since $\epsilon_n^{(j)}>0$ for all
$j$. Then as a direct consequence of (\ref{eq:approx}), we have for
$nt\in \nat$,
\begin{eqnarray}
\label{eq:approx2}\A^{(i)}(t) + t\eta(\delta) = \Nbf^{(i)}(t)
+t\delta + \sum_{j=1}^K
\Pbf^{(j,i)}_{n,\delta}(\D^{(j)}(t)+t\eta(\delta)).
\end{eqnarray}

If $\Nbf_{n,\delta}$ is the polygonal approximation of
$t\rightarrow\Nbf(t)+t\delta$ with step $1/n$, we have clearly
$\dot{\Nbf}_{n,\delta}\rightarrow \dot{\Nbf}+\delta$ as $n$ tends to
infinity. Similarly, we have as $n$ tends to infinity,
\begin{eqnarray*}
\dot{\Pbf}_{n,\delta}^{(j,i)}(\D^{(j)}(t)+t\eta(\delta))\rightarrow
\left\{\begin{array}{ll}
(1-\epsilon^{(j)}(t))\dot{\Pbf}^{(j,i)}(\D^{(j)}(t))\frac{(\dot{\D}^{(j)}(t)+\eta(\delta))}{\dot{\D}^{(j)}(t)}&\\
\quad \quad
+\epsilon^{(j)}(t)R^{(j,i)}(\dot{\D}^{(j)}(t)+\eta(\delta))&\mbox{if
$\dot{\D}^{(j)}(t)>0$},\\
R^{(j,i)}\eta(\delta)&\mbox{otherwise,}
\end{array} \right.
\end{eqnarray*}
where $\epsilon^{(j)}(t) =
\eta^{(j)}(\delta)/(\eta^{(j)}(\delta)+\dot{\D}^{(j)}(t))<1$. Hence
when $n$ tends to infinity and $\delta$ tends to zero, we have
$\dot{\Nbf_{n,\delta}}\rightarrow \dot{\Nbf}$ and
$\dot{\Pbf}_{n,\delta}^{(j,i)}\rightarrow\dot{\Pbf}^{(j,i)}$.

We consider now the queueing equation $\D=\Phi(\A,\Sbf)$ and
construct the approximating sequence for $\Sbf$.

We begin with a first general lemma: given three processes $\A \leq
\D$ and $\Sbf$, we construct a piecewise linear function $\Sbf_n$
(with step $1/n$) as follows (with $nt\in \nat$):
\begin{itemize}
\item if $\A(t)=\D(t)$ and $\A(t+1/n)=\D(t+1/n)$, then
$\Sbf_n(t+1/n)-\Sbf_n(t)=\Sbf(t+1/n)-\Sbf(t)$;
\item otherwise, $\Sbf_n(t+1/n)-\Sbf_n(t)=\D(t+1/n)-\D(t)$.
\end{itemize}
We will denote this construction by $\Sbf_n = \Upsilon_n(\A,
\D,\Sbf)$.
\begin{lem}\label{lem:approx}
Let $(\A,\D,\Sbf)$ be absolutely continuous functions of
$\Cbb(\real_+^K)\times \Cbb_0(\real_+^K) \times
\Cbb_0(\real_+^K)$such that $\Phi(\A,\Sbf)=\D$. %Assume that there exists a sequence such that $\dot{\A}_n\rightarrow \dot{\A}$, $\dot{\D}_n\rightarrow \dot{\D}$ and $\A_n\geq \D_n$ for all $n$.
We denote $\Sbf_n= \Upsilon_n(\A, \D,\Sbf)$. We have
$\D_n=\Phi(\A_n,\Sbf_n)$ where $(\A_n,\D_n)$ is the polygonal
approximation of $(\A,\D)$ with step $1/n$ and we have the following
convergence as $n$ tends to infinity: $\Sbf_n\to \Sbf$,
$\dot{\Sbf}_n\rightarrow \dot{\Sbf}$ and $\int_0^\infty
\Isf^{\Sbf}(\dot{\Sbf}_{n}(t))dt\to\int_0^\infty
\Isf^{\Sbf}(\dot{\Sbf}(t))dt$ .
\end{lem}
\proof We denote $\tilde{\D}_n=\Phi(\A_n,\Sbf_n)$. From the proof of
Lemma \ref{lem:trans}, we have
\begin{eqnarray*}
\tilde{\D}_n(t+1/n)-\tilde{\D}_n(t) = \inf_{t\leq s\leq
t+1/n}\left\{ \Sbf_n(t+1/n)
-\Sbf_n(s)+\A_n(s)-\tilde{\D}_n(t)\right\}\wedge
\left\{\Sbf_n(t+1/n)-\Sbf_n(t) \right\},
\end{eqnarray*}
since all the functions are linear on the interval $(t,t+1/n)$, we
have (with $nt\in \nat$),
\begin{eqnarray*}
\tilde{\D}_n(t+1/n)-\tilde{\D}_n(t) &=& \left\{
\A_n(t+1/n)-\tilde{\D}_n(t)\right\}\wedge
\left\{\Sbf_n(t+1/n)-\Sbf_n(t) \right\}.
\end{eqnarray*}
If $\tilde{\D}_n(t)=\D_n(t)$, then we have clearly
$\tilde{\D}_n(t+1/n)=\D_n(t+1/n)$ since
\begin{itemize}
\item if $\A_n(t)=\D_n(t)$ and $\A_n(t+1/n)=\D_n(t+1/n)$, then we
have $\Sbf(t+1/n)-\Sbf(t)\geq
\D_n(t+1/n)-\D_n(t)=\A_n(t+1/n)-\tilde{\D}_n(t)$ see
(\ref{ineq:lem}) for the inequality;
\item otherwise, $\Sbf_n(t+1/n)-\Sbf_n(t)=\D_n(t+1/n)-\D_n(t)$ by definition and $\A_n(t+1/n)\geq
\D_n(t+1/n)$.
\end{itemize}
This proves the first part of the lemma. Moreover it follows
directly form the definition of $\Upsilon$ that
$\Sbf_n\left(t+\frac{1}{n}\right)-\Sbf_n(t) \leq
\Sbf\left(t+\frac{1}{n}\right)-\Sbf(t)$, hence we have for all $t$,
$\limsup_{n\to \infty}\Sbf_n(t)\leq \Sbf(t)$ by a continuity
argument. The fact that $\Sbf_n\to\Sbf$ follows directly from
Fatou's Lemma and the fact that $\dot{\Sbf}_n\to\dot{\Sbf}$. We now
prove this last fact, let $C=\{t,\: \A(t)=\D(t)\}$. $C$ is a closed
set and according to Lemma \ref{lem:propPhi}, we have for all $t\in
C^c$ (the complementary set of $C$), $\dot{\Sbf}(t)=\dot{\D}(t)$.
For such $t\in C^c$, we have for $\epsilon>0$ sufficiently small and
for sufficiently large $n$, $\A_n(u)\neq \D_n(u)$ for all $|u-t|\leq
\epsilon$. Hence we have $\dot{\Sbf}_n(t)=\dot{\D}_n(t)\rightarrow
\dot{\D}(t)$. Now for $t\in C^o$ in the interior of $C$, we have
clearly $\dot{\Sbf}_n(t)\rightarrow \dot{\Sbf}(t)$. Hence we have
$\dot{\Sbf}_n(t)\rightarrow \dot{\Sbf}(t)$ for $t\in C^o\cup C^c$.

We prove the last statement of the lemma. Since any open set of
$\real$ is a countable union of disjoint intervals,
\begin{eqnarray*}
\int_{C^c} \Isf^{\Sbf}(\dot{\Sbf}_n(t))dt = \int_{C^c}
\Isf^{\Sbf}(\dot{\D}_n(t))dt &\leq& \int_{C^c}
\Isf^{\Sbf}(\dot{\D}(t))dt,\quad \mbox{by Jensen's
inequality}\\
&=&\int_{C^c} \Isf^{\Sbf}(\dot{\Sbf}(t))dt,
\end{eqnarray*}
and also directly still by Jensen's inequality $\int_{C^o}
\Isf^{\Sbf}(\dot{\Sbf}_n(t))dt\leq \int_{C^o}
\Isf^{\Sbf}(\dot{\Sbf}(t))dt$. The convergence then follows from
\begin{eqnarray*}
\liminf_{n\rightarrow \infty}\int_0^\infty
\Isf^{\Sbf}(\dot{\Sbf}_{n}(t))dt&\geq&
\int_0^\infty\liminf_{n\rightarrow
\infty} \Isf^{\Sbf}(\dot{\Sbf}_{n}(t))dt\\
&\geq&\int_0^\infty \Isf^{\Sbf}(\dot{\Sbf}(t))dt,
\end{eqnarray*}
where the first inequality is due to Fatou's Lemma and the second
one to the lower semicontinuity of $\Isf^{\Sbf}$.
\endproof

We define the sequence
$\JN_{n,\delta}=(\Sbf_{n,\delta},\Pbf_{n,\delta},\Nbf_{n,\delta})$
where $\Sbf_{n,\delta}(t)=\Upsilon_n(\A(t)+\eta t, \D(t)+\eta
t,\Sbf(t)+\eta t)$. Note that we have $\D(t) +\eta t=\Phi(\A(t)+\eta
t,\Sbf(t) +\eta t)$, hence Lemma \ref{lem:approx} applies, in
particular, we have $\dot{\Sbf}_{n,\delta}(t)\rightarrow
\dot{\Sbf}(t)+\eta(\delta)$ as $n$ tends to infinity.

We have $\JN_{n,\delta}\in \Dcal_\JN$ by construction and the
sequence $\{\JN_{n,\delta_n}\}_n$ satisfies (\ref{eq:S2}) for some
$\delta_n\rightarrow 0$. Moreover, we have thanks to
(\ref{eq:approx2}) and Lemma \ref{lem:approx},
\begin{eqnarray*}
\left\{\begin{array}{lcl}\A_{n,\delta}&=&\Gamma(\D_{n,\delta},\JN_{n,\delta}),\\
\D_{n,\delta} &=&
\Phi(\A_{n,\delta},\JN_{n,\delta}),\end{array}\right.\Leftrightarrow
(\A_{n,\delta},\D_{n,\delta})=\Psi(\JN_{n,\delta}),
\end{eqnarray*}
where $\A_{n,\delta}$ and $\D_{n,\delta}$ are the polygonal
approximation of $\A(t)+\eta t$ and $\D(t)+\eta t$ with step $1/n$
and $\Psi$ has been defined in Section \ref{sec:extpsi}.

For $n\rightarrow \infty$ and $\delta\rightarrow 0$, we have
$(\A_{n,\delta},\D_{n,\delta})\rightarrow (\A,\D)$, hence we have
$\Psi(\JN_{n,\delta})\rightarrow (\A,\D)$, i.e. the sequence
$\{\JN_{n,\delta_n}\}_n$ satisfies (\ref{eq:S3}).

We now show that (\ref{eq:S4}) is also satisfied. We fix $T>0$ and
prove first that we have, for $\delta$ sufficiently small,
\begin{eqnarray}
\label{eq:approxI}\lefteqn{\int_0^T
\Isf^{\Sbf}(\dot{\Sbf}(t))+\tilde{\Dsf}(\dot{\Pbf}(t)\|R)+\Isf^{\Nbf}(\dot{\Nbf}(t))dt
-er(\delta)T}\\
\nonumber&\leq& \liminf_{n\to \infty} \int_0^T
\Isf^{\Sbf}(\dot{\Sbf}_{n,\delta}(t))+\tilde{\Dsf}(\dot{\Pbf}_{n,\delta}(t)\|R)+\Isf^{\Nbf}(\dot{\Nbf}_{n,\delta}(t))dt\\
\nonumber&\leq& \limsup_{n\to \infty} \int_0^T
\Isf^{\Sbf}(\dot{\Sbf}_{n,\delta}(t))+\tilde{\Dsf}(\dot{\Pbf}_{n,\delta}(t)\|R)+\Isf^{\Nbf}(\dot{\Nbf}_{n,\delta}(t))dt\\
\label{eq:approxII}&\leq& \int_0^T
\Isf^{\Sbf}(\dot{\Sbf}(t))+\tilde{\Dsf}(\dot{\Pbf}(t)\|R)+\Isf^{\Nbf}(\dot{\Nbf}(t))dt
+er(\delta)T,
\end{eqnarray}
where $er(\delta)$ tends to zero as $\delta$ tends to zero, from
which (\ref{eq:S4}) follows by monotonicity.

We first deal with the case of the sequence of processes
$\{\Sbf_{n,\delta}\}_n$ (we can restrict ourselves to the one
dimensional case). We denote $\Sbf_\delta(t)=\Sbf(t)+\eta(\delta)t$.

%We first take $T>0$.
We define $\varsigma=\esssup\{\dot{\Sbf}(t),\: t\leq
T\}=\inf\{u,\:Leb[t\leq T,\:\dot{\Sbf}(t)>u ]=0\}$, where $Leb$ is
for the Lebesgue measure. Since $\int_0^T
\Isf^\Sbf(\dot{\Sbf}(t))dt<\infty$, $\varsigma$ belongs to the
domain of $\Isf^\Sbf$ which is open on the right. Hence we can find
$\epsilon>0$ such that $\varsigma+\epsilon$ still belongs to this
domain and take $\delta$ such that $\eta(\delta)<\epsilon$.
Moreover, since $\Isf^\Sbf$ is convex, it is uniformly continuous on
$[0,\varsigma+\epsilon]$. Hence, we can assume that we have
$\beta(\alpha)\to 0$ as $\alpha\to 0$ such that,
\begin{eqnarray*}
\forall x,y\in [0,\varsigma+\epsilon], \:|x-y|<\alpha \Rightarrow
|\Isf^\Sbf(x)-\Isf^\Sbf(y)|\leq \beta(\alpha).
\end{eqnarray*}

From Lemma \ref{lem:approx}, we have
\begin{eqnarray}
\nonumber\int_0^T \Isf^{\Sbf}(\dot{\Sbf}(t))dt -\beta(\eta(\delta)) T\leq&&\\
\nonumber\lim_{n\rightarrow \infty}\int_0^T
\Isf^{\Sbf}(\dot{\Sbf}_{n,\delta}(t))dt&=& \int_0^T
\Isf^{\Sbf}(\dot{\Sbf}_\delta(t))dt\\
\label{eq:ineqsup}&&\leq \int_0^T \Isf^{\Sbf}(\dot{\Sbf}(t))dt
+\beta(\eta(\delta)) T.
\end{eqnarray}
\iffalse
and the converse inequality,
\begin{eqnarray*}
\liminf_{n\rightarrow \infty}\int_0^\infty
\Isf^{\Sbf}(\dot{\Sbf}_{n,\delta}(t))dt &\geq&\int_0^T
\Isf^{\Sbf}(\dot{\Sbf}_\delta(t))dt\geq \int_0^T
\Isf^{\Sbf}(\dot{\Sbf}(t))dt -\beta(\eta(\delta)) T.
\end{eqnarray*}
\fi Hence we proved (\ref{eq:approxI}) and (\ref{eq:approxII}) for
$\Isf^\Sbf$.

The same kind of arguments can be repeated for $\Nbf_{n,\delta}$
which is just the polygonal approximation of $t\mapsto\Nbf(t)+\delta
t$. Note that $\{\JN_{n,\delta}\}_n\in\Dcal_\JN^\nat$ implies that
$\Nbf^{(i)}_{n,\delta}(t)=0$ for all $i\notin \Scal$. For $i\in
\Scal$, we can use the fact that the domain of $\Isf^{\Nbf^{(i)}}$
is open as previously. In the case of $\Pbf_{n,\delta}$, we can not
use the argument on the openness of the domain, but we have
$\tilde{\Dsf}(R^{(i)}\|R^{(i)})=0$ and then the convexity of
$\tilde{\Dsf}$ directly implies that
$\tilde{\Dsf}(\dot{\Pbf}_{n,\delta}^{(i)}\|R^{(i)})\leq
\tilde{\Dsf}(\dot{\Pbf}^{(i)}\|R^{(i)})$, from which we derive an
equivalent of (\ref{eq:ineqsup}).
%\endproof

\subsection{Exponential tightness}

We first recall some definitions. A sequence of random variables
$\{X_n\}_n\in (\real^K)^\nat$ is exponentially tight if
\begin{eqnarray*}
\lim_{M\rightarrow \infty}\limsup_{n\rightarrow \infty} \frac{1}{n}
\log \prob(\|X_n\|>M)= -\infty.
\end{eqnarray*}
For $\delta>0$ and $T>0$, define the modulus of continuity in
$\Dbb(E)$ by
\begin{eqnarray*}
w'(\Xbf,\delta,T):= \inf_{\{t_i\}}\max_i
\sup_{s,t\in[t_{i-1},t_i)}d(\Xbf(s),\Xbf(t)),
\end{eqnarray*} where
the infimum is over $\{t_i\}$ satisfying
\begin{eqnarray*}
0=t_0< t_1<\dots <t_{m-1}<T\leq t_m
\end{eqnarray*}
and $\min_{1\leq i\leq n}(t_i-t_{i-1})>\delta$.

Theorem 4.1 of \cite{fenkur} tells us: let $\Tcal_0$ be a dense
subset of $\real_+$. Suppose that for each $t\in \Tcal_0$,
$\{\Xbf_n(t)\}_n$ is exponentially tight. Then $\{\Xbf_n\}_n$ is
exponentially tight in $\Dbb(E)$ if and only if for each
$\epsilon>0$ and $T>0$,
\begin{eqnarray}
\label{exptight}\lim_{\delta\rightarrow 0}\limsup_{n\rightarrow
\infty}\frac{1}{n}\log \prob(w'(\Xbf_n,\delta,T)>\epsilon)=-\infty.
\end{eqnarray}

A sequence of stochastic processes $\{\Xbf_n\}_n$ that is
exponentially tight in $\Dbb(E)$ is $C$-exponentially tight if for
each $\eta>0$ and $T>0$,
\begin{eqnarray}
\label{Cexptight}\limsup_{n\rightarrow\infty} \frac{1}{n}
\log\prob(\sup_{s\leq T} d(\Xbf_n(s),\Xbf_n(s-))\geq \eta)=-\infty.
\end{eqnarray}

Then Theorem 4.13 of \cite{fenkur} gives: an exponentially tight
sequence $\{\Xbf_n\}_n$ in $\Dbb(E)$ is $C$-exponentially tight if
and only if each rate function $\Isf$ that gives the LDP for a
subsequence $\{\Xbf_{n(k)}\}_{n(k)}$, satisfies $\Isf(\x)=\infty$
for each $\x\in \Dbb(E)$ such that $\x\notin \Cbb(E)$.

The stochastic assumptions of Section \ref{sec:sto} ensure that the
sequence of processes $\{\JN_n\}_n$ satisfies a LDP with good rate
function (this implies that the sequence is exponentially tight)
giving an infinite mass to discontinuous path. Hence the sequence of
processes $\{\JN_n\}_n$ is $C$-exponentially tight.

We have to show that the sequence of processes $\{(\A_n,\D_n)\}_n$
is exponentially tight. The fact of dealing with non-decreasing
processes simplifies the definitions. For $\Xbf\in \Dbb(\real_+^K)$
(or $\Dbb(\Mbb^K)$) non-decreasing, $\delta>0$ and $T>0$, we define
$w_\delta(\Xbf,T)=\sup_{t\in[0,T]}\|\Xbf(t+\delta)-\Xbf(t)\|$. We
have clearly $w'(\Xbf,\delta,T) = w_\delta(\Xbf,T)$ and if
$\{\Xbf_n(0)\}_n$ is exponentially tight then (\ref{exptight})
implies that $\{\Xbf_n(t)\}_n$ is exponentially tight for each
$t>0$. Lemmas \ref{lem:wPhi} and \ref{lem:wGamma} show that
conditions (\ref{exptight}) and (\ref{Cexptight}) are satisfied for
the sequence of processes $\{(\A_n,\D_n)\}_n$. The exponential
tightness of $\{(\A_n(0),\D_n(0))\}_n$ is clear since
$\A_n(0)=\D_n(0)=0$.
%follows directly form $\D_n(0)=0$ and $\A_n(0)=\Nbf_n(0)$ and the fact that $\{\Nbf_n(0)\}_n$ satisfies a LDP with rate function $I^0$.

\subsection{Large deviations results}

\begin{prop}\label{prop:sldp}
The sequence of processes $\{(\A_n, \D_n)\}_n$ satisfies a LDP in
$\Dbb(\real^K_+)\times \Dbb(\real_+^K)$ with good rate function
$\Isf^{\A,\D}$. For $\A,\D$ absolutely continuous and such that
$\A(0)=\D(0)=0$ and $\A\geq \D$, $\Isf^{\A,\D}$ is given by
\begin{eqnarray}
\label{eq:LDP}\Isf^{\A,\D}(\A, \D)=\int_0^\infty
H(\A(s),\D(s),\dot{\A}(s),\dot{\D}(s))ds,
\end{eqnarray}
where $H(A,D,\dot{A},\dot{D}):=\inf_{P,N}
h(A,D,\dot{A},\dot{D},P,N)$, with $h$ given by,
\begin{eqnarray*}
\lefteqn{h(A,D,\dot{A},\dot{D},P,N):=}\\
&&\sum_{i\in
E(A,D)}\Isf^{\Sbf^{(i)}}(\dot{D}^{(i)})\ind_{\{\dot{D}^{(i)}>\mu^{(i)}\}}+\sum_{i\notin
E(A,D)}\Isf^{\Sbf^{(i)}}(\dot{D}^{(i)}) +\sum_i \dot{D}^{(i)}
\tilde{\Dsf}(P^{(i)}\|R^{(i)})+\Isf^\Nbf(N)
\end{eqnarray*}
where $E(A,D)=\{i,\: A^{(i)}=D^{(i)}\}$ and with the infimum taken
over the set of $(P,N)\in \Mbb^K\times \real_+^K$ such that
\begin{eqnarray*}
\dot{A} &=& N+P^t\dot{D}.%\quad \mbox{and,}\quad \rho(P)<1.
\end{eqnarray*}
For all other $\A,\D$, we have $\Isf^{\A,\D}(\A,\D)=\infty$.
\end{prop}
\proof Thanks to the results of previous sections, conditions of
Proposition \ref{prop:contraction} are satisfied and we define
\begin{eqnarray}
\label{eq:ratetilde}\tilde{\Isf}^{\A,\D}(\A, \D)=\inf
\left\{\lim_{n\rightarrow \infty}\Isf^\JN(\JN_n),\: \{\JN_n\}_n\in
\Scal(\A,\D)\right\},
\end{eqnarray}
where we recall that $\Scal(\A,\D)=\cup_\JN \Scal(\JN,\A,\D)$, and
$\Scal(\JN,\A,\D)$ is defined in Proposition \ref{prop:JNapprox}. We
have to show that $\tilde{\Isf}^{\A,\D}=\Isf^{\A,\D}$ given by
(\ref{eq:LDP}).

\iffalse We define
\begin{eqnarray}
\label{eq:ratetilde}\tilde{I}^{\A,\D}(\A, \D)=\inf
\left\{I^\JN(\JN),\: (\A,\D) =
(\Gamma(\D,\JN),\Phi(\A,\JN))\right\}.
\end{eqnarray}
We have to show that $\tilde{I}^{\A,\D}=I^{\A,\D}$ given by
(\ref{eq:LDP}). Note that $\tilde{I}^{\A,\D}$ is the lower
semicontinuous regularization of $\bar{I}^{\A,\D}$ defined as in
(\ref{eq:ratetilde}) but where the infimum is restricted to
piece-wise linear Jackson networks with routing matrix having a
derivative of spectral radius less than one. Hence it is sufficient
to prove that $I^{\A,\D}=\bar{I}^{\A,\D}$ and to show that
$I^{\A,\D}$ is a good rate function. \fi

Consider $\JN\in \Dcal_\JN$ and let $(\A,\D)=\Psi(\JN)$. Let
$\tau=\{0=t_0<t_1<\dots \}$ be such that the processes $\A,\D, \Sbf,
\Nbf$ and $\D\circ\Pbf$ have a constant derivative on each
$(t_k,t_{k+1})$. %We denote the corresponding derivative on this interval by $\dot{\A}_k, \dot{\D}_k\dots$
%such that $I^\JN(\JN)<\infty$ and such that there exists a solution $(\A,\D)$ to the fixed point equation (\ref{fixpointeq}).
Then from $\A=\Gamma(\D,\JN)$, we derive
\begin{eqnarray*}
\dot{\A}^{(i)}(t) = \dot{\Nbf}^{(i)}(t) +\sum_j \dot{\D}^{(j)}(t)
\dot{\Pbf}^{(j,i)}(\D^{(j)}(t)).
\end{eqnarray*}
From $\D = \Phi(\A,\Sbf)$, we get the following constraints:
\begin{itemize}
\item if $\A^{(i)}(t_k)>\D^{(i)}(t_k)$ or $\A^{(i)}(t_{k+1})>\D^{(i)}(t_{k+1})$, then we have
$\dot{\D}^{(i)}(t) = \dot{\Sbf}^{(i)}(t)$ for $t\in (t_k,t_{k+1})$;
\item otherwise $\A^{(i)}(t)=\D^{(i)}(t)$ for $t\in (t_k,t_{k+1})$ and we have $\dot{\Sbf}^{(i)}(t)\geq
\dot{\A}^{(i)}(t)=\dot{\D}^{(i)}(t)$ for $t\in (t_k,t_{k+1})$.
\end{itemize}

Now we can compute $\Isf^\JN(\JN)$ as follows
\begin{eqnarray*}
\Isf^\JN(\JN) %&=& \sum_k \left\{I^{\Sbf}(\dot{\Sbf}_k)+\sum_i D(\dot{\Pbf}^{(i)}(\D^{(i)}_k)\|R^{(k)})+I^\Nbf(\dot{\Nbf}_k)\right\}(t_{k+1}-t_k)\\
%&=&\sum_k\left\{I^{\Sbf}(\dot{\Sbf}_k)+I^\Nbf(\dot{\Nbf}_k)\right\}(t_{k+1}-t_k)+\sum_k\sum_iD(\dot{P}^{(i)}(D^{(i)}(t_k))\|R^{(k)})(t_{k+1}-t_k)\\
&=&\int_0^\infty\sum_{i\in
E(A,D)}\Isf^{\Sbf^{(i)}}(\dot{\Sbf}^{(i)}(s))+\sum_{i\notin
E(A,D)}\Isf^{\Sbf^{(i)}}(\dot{\D}^{(i)}(s))+\Isf^\Nbf(\dot{\Nbf}(s)) ds\\
&&+\int_0^\infty\sum_j \dot{\D}^{(j)}(s)
\tilde{\Dsf}(\dot{\Pbf}^{(j)}(s)\|R^{(j)})ds\\
&\geq& \int_0^\infty
h(\A(s),\D(s),\dot{\A}(s),\dot{\D}(s),\dot{\Pbf}(s),\dot{\Nbf}(s))ds
\geq \Isf^{\A,\D}(\Psi(\JN)),
\end{eqnarray*}
since for $i\in E(\A(s),\D(s))$, we have
$\Isf^{\Sbf^{(i)}}(\dot{\Sbf}^{(i)}(s))\geq
\Isf^{\Sbf^{(i)}}(\dot{\D}^{(i)}(s))\ind_{\{\dot{\D}^{(i)}(s)>\mu^{(i)}\}}$
because $\dot{\Sbf}^{(i)}(s)\geq \dot{\D}^{(i)}(s)$ and
$\Isf^{\Sbf^{(i)}}$ is non-negative, convex with $\mu^{(i)}$ as
unique zero. Hence, we have $\tilde{\Isf}^{\A,\D}\geq \Isf^{\A,\D}$.

Consider now $(\A,\D)$ such that $\Isf^{\A,\D}(\A,\D)<\infty$, %and let $(\A_n,\D_n)$ be the respective polygonal approximations with step $1/n$. Then  we have $I^{\A,\D}(\A_n,\D_n)\rightarrow I^{\A,\D}(\A,\D)$.
then we denote by $(\pbf(s),\nbf(s))$ the argument that achieves the
minimum in $H(\A(s),\D(s),\dot{\A}(s),\dot{\D}(s))$ for any fixed
$s$ (note that $h$ is a good rate function). Let $\Pbf(\D(t))
=\int_0^t \pbf(s) ds$ and $\Nbf(t) = \int_0^t \nbf(s) ds$, note that
$\pbf$ and $\nbf$ are measurable since $H$ is a good rate function.
We have $\A = \Gamma(\D,\Pbf,\Nbf)$. Now define $\sbf(s)$ as
follows:
\begin{itemize}
\item if $\A^{(i)}(s) =\D^{(i)}(s)$ then $\sbf^{(i)}(s)=\dot{\D}^{(i)}(s)
\vee \mu^{(i)}$;
\item if $\A^{(i)}(s) > \D^{(i)}(s)$ then
$\sbf^{(i)}(s)=\dot{\D}^{(i)}(s)$.
\end{itemize}
We have $\D=\Phi(\A,\Sbf)$ with $\Sbf(t)=\int_0^t \sbf(s) ds$. Hence
we have $(\A,\D) = (\Gamma(\D,\JN),\Phi(\A,\JN))$ for
$\JN=(\Sbf,\Pbf,\Nbf)$ and
$\Isf^\JN(\JN)=\Isf^{\A,\D}(\A,\D)<\infty$ by construction. Hence
the sequence $\Scal(\JN,\A,\D)=\{\JN_n\}_n$ is well-defined and we
have $\tilde{\Isf}^{\A,\D}(\A,\D)\leq
\lim_{n\rightarrow\infty}\Isf^\JN(\JN_n)=\Isf^\JN(\JN)=\Isf^{\A,\D}(\A,\D)$.
%,where the last fact follows directly from the construction of $\JN$.
\endproof

From this proposition, it is quite easy to derive a LDP for the
process $\Q_n(t) :=\A_n(t)-\D_n(t)$ counting the number of customers
in each queue. Thanks to the contraction principle, we have
\begin{eqnarray*}
\Isf^\Q(\Q) = \inf\{\Isf^{\A,\D}(\A,\D), \Q=\A-\D \},
\end{eqnarray*}
which gives directly Theorem \ref{th:qu}.

\subsection*{Acknowledgements}
I am thankful to Anatolii Puhalskii for insightful comments and for
providing me with a copy of \cite{Puh:jack}

\section{Appendix}

\subsection{Properties of the map $\Gamma$ and $\Phi$}
For $\Xbf\in \Dbb(\real_+^K)$, $\delta>0$ and $T>0$, we define
$w_\delta(\Xbf,T)=\sup_{t\in[0,T]}\|\Xbf(t+\delta)-\Xbf(t)\|$.
\begin{lem}\label{lem:wPhi}
We have
\begin{eqnarray*}
w_\delta(\Phi(\Xbf,\Ybf),T)\leq w_\delta(\Ybf,T).
\end{eqnarray*}
\end{lem}
\proof It is clearly sufficient to consider the case $K=1$. We will
prove that
\begin{eqnarray}
\label{ineq:lem}\Phi(\Xbf,\Ybf)(t+\delta)-\Phi(\Xbf,\Ybf)(t) \leq
\Ybf(t+\delta)-\Ybf(t),
\end{eqnarray}
from which the lemma follows. If $\Phi(\Xbf, \Ybf)(t)=\Ybf(t)$, then
we have $\Phi(\Xbf, \Ybf)(t+\delta) \leq \Ybf(t+\delta)$ and
(\ref{ineq:lem}) is clear.

Assume now that $\Phi(\Xbf, \Ybf)(t) =\inf_{0\leq
s<t}\left\{\Ybf(t)-\Ybf(s)+\Xbf(s)\right\}<\Ybf(t)$. %In particular, we have $t-s+\Ybf^{\leftarrow}(\Xbf(s))\leq \Ybf^{\leftarrow}(\Xbf(t))$ for all $s\leq t$. Moreover we have
We have
\begin{eqnarray*}
\Ybf(t+\delta)-\Ybf(s)+\Xbf(s)=
\Ybf(t)-\Ybf(s)+\Xbf(s)+\Ybf(t+\delta)-\Ybf(t),
\end{eqnarray*}
and (\ref{ineq:lem}) follows by taking the minimum in $s\in[0,t]$
and observing that $\Phi(\Xbf, \Ybf)(t+\delta)\leq \inf_{0\leq
s<t}\left\{\Ybf(t+\delta)-\Ybf(s)+\Xbf(s)\right\}$.
\endproof

The following lemma is clear:
\begin{lem}\label{lem:wGamma}
We have
\begin{eqnarray*}
w_\delta(\Gamma(\Xbf,\Pbf,\Nbf),T)\leq
w_\delta(\Nbf,T)+w_\delta(\Pbf,\|\Xbf(T)\|).
\end{eqnarray*}
\end{lem}

\begin{lem}\label{lem:propPhi}
Assume $\Sbf\in \Dbb_0(\real_+)$ is absolutely continuous, then for
any $\A\in \Dbb(\real_+)$, we have $\D:=\Phi(\A, \Sbf)$ is
absolutely continuous and,
\begin{itemize}
\item for all $t$ such that
$\A(t)>\D(t)$, we have $\dot{\D}(t)=\dot{\Sbf}(t)$;
\item if $\A(t)=\D(t)$ for $t\in (u,v)$ with $u<v$, then we have $\dot{\Sbf}(t)\geq
\dot{\A}(t)=\dot{\D}(t)$ for $t\in (u,v)$.
\end{itemize}
\end{lem}
\proof It follows directly form (\ref{ineq:lem}) that if $\Sbf$ is
absolutely continuous, then $\Phi(\Xbf,\Sbf)$ is absolutely
continuous for any $\Xbf$. The rest of the lemma is obvious.
\endproof

\subsection{Auxiliary results}
\begin{lem}\label{lem:specrad}
Given a substochastic matrix $R$ such that $\rho(R)<1$ and a
substochastic matrix $P$ such that the support of $P$ is included in
the support of $R$, i.e. $R^{(i,j)}=0\Rightarrow P^{(i,j)}=0$. Then
for any $\epsilon$ such that $0< \epsilon^{(i)}\leq 1$ for all $i$,
the matrix with coefficients $M^{(i,j)}=(1-\epsilon^{(i)})
P^{(i,j)}+\epsilon^{(i)} R^{(i,j)}$ is of spectral radius less than
$1$.
\end{lem}
\proof By a suitable permutation of rows and columns, we can assume
that $R$ is given in its canonical form
\begin{eqnarray}
\label{eq:canonic}R =\left(\begin{array}{ccccc}
S_1(R)&*&*&*\\
0&S_2(R)&*&*\\
0&0&\ddots&*\\
0&0&0&S_n(R)\end{array} \right),
\end{eqnarray}
where each $S_i(R)$ is an irreducible matrix. We have $\rho(R)<1$ if
and only if each $S_i(R)$ is not a stochastic matrix.

In view of the assumption on the support of $P$, the matrix $P$ has
the same structure as (\ref{eq:canonic}) and we have with the same
notation as above, $S_i(M)$ which is an irreducible and not
stochastic matrix.
\endproof

\subsection{An example}
In this section, we construct 2 different sequences of Jackson
networks $\JN_n^1$ and $\JN_n^2$ such that their fluid limits are
the same
\begin{eqnarray*}
\JN_n^1\rightarrow \JN \quad \mbox{and}\quad \JN_n^2\rightarrow \JN,
\end{eqnarray*}
but such that
\begin{eqnarray*}
(\A_n^1,\D_n^1)=\Psi(\JN_n^1)&\rightarrow& (\A^1,\D^1),\\
(\A_n^2,\D_n^2)=\Psi(\JN_n^2)&\rightarrow& (\A^2,\D^2),
\end{eqnarray*}
with $(\A^1,\D^1)\not= (\A^2,\D^2)$.

\begin{figure}[hbt]
\begin{picture}(0,120)
\psfrag{sigma0}{$\Nbf$} \psfrag{sigma1}{$\Sbf$}
\psfrag{p11}{$\Pbf^{(1,1)}$} \psfrag{p12}{$\Pbf^{(1,2)}$}
\put(-200,0){\epsfig{file=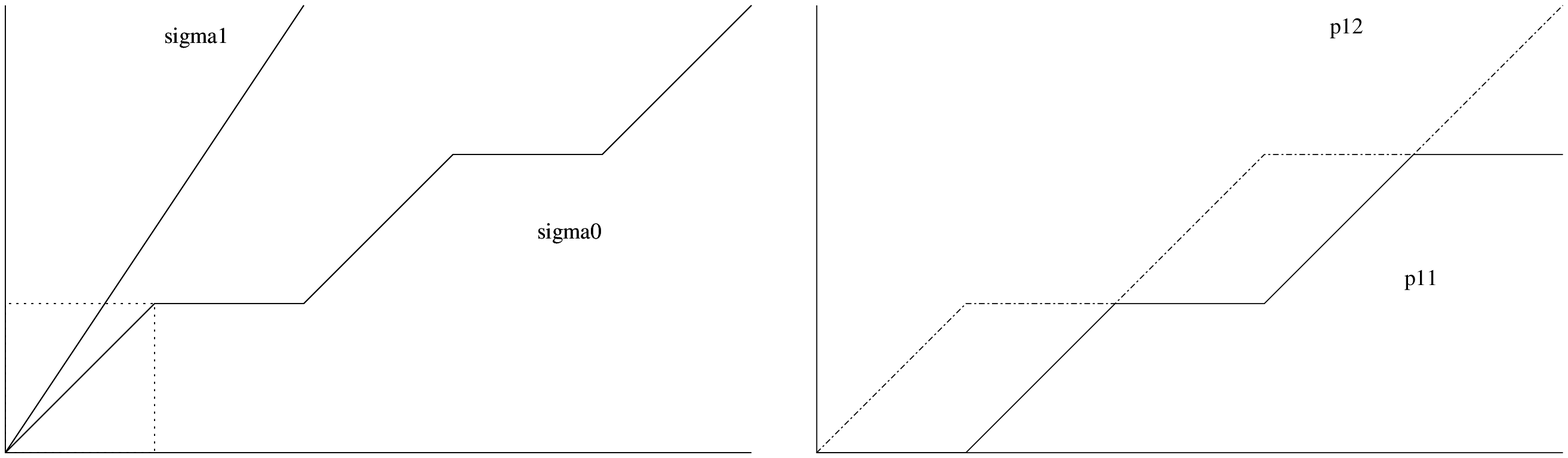 ,width =14cm}}
\end{picture}
\caption{Fluid networks: $\JN$} \label{fig:JN}
\end{figure}

\begin{figure}[hbt]
\begin{picture}(0,140)
\psfrag{Q}{$\Q^2$} \psfrag{D2}{$\D^2$} \psfrag{x}{$x\alpha$}
\put(-100,0){\epsfig{file=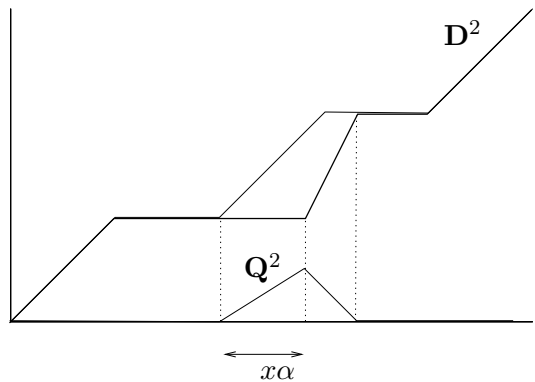 ,width =7cm}}
\end{picture}
\caption{Departure process and queue length process} \label{fig:dep}
\end{figure}

We consider a toy example with only one station (hence we omit the
superscript $.^{(1)}$ that refers to that only station). Once a
customer is served, he can either go out of the network or go back
to this same node. We define the following driving sequences:
\begin{eqnarray*}
T^{n} &=&
(\underbrace{1,\dots,1}_{n},n,\underbrace{1,\dots,1}_{n},n,\dots),\\
\sigma^{n}&=& \alpha(1,1,\dots),
\end{eqnarray*}
with $\alpha<1$. We define now two different routing sequences
\begin{eqnarray*}
\nu^{n} &=&
(\underbrace{2,\dots,2}_{n+1},\underbrace{1,\dots,1}_{n+1},\dots),\\
\nu^{n}(x) &=& (\underbrace{2,\dots,2}_{\lfloor
  xn\rfloor},1,\underbrace{2,\dots,2}_{n-\lfloor xn
  \rfloor},\underbrace{1,\dots,1}_{\lfloor
  xn\rfloor},2,\underbrace{1,\dots,1}_{n-\lfloor xn\rfloor}, \dots),
\end{eqnarray*}
where $x<1$. We denote by $\JN_n^1=\left\{ \sigma^n,\nu^n,T^n
\right\}$ and $\JN_n^2=\left\{ \sigma^n,\nu^n(x),T^n \right\}$.
$\nu^n(x)$ is obtained from $\nu^n$ by only interchanging a 1 and a
2. Hence we have
\begin{eqnarray*}
\JN_n^1\rightarrow \JN \quad \mbox{and}\quad \JN_n^2\rightarrow \JN.
\end{eqnarray*}
Indeed the fluid network $\JN$ is given on Figure \ref{fig:JN}.
\\
In the fluid limit, in case 1, the queue is always empty and the
departure process is the same as the arrival process $\Nbf$. In case
2, the fluid limit of the departure process and the queue length
process is given on Figure \ref{fig:dep}.
\\
To explain $\D^2$, we write for each arrival (number on the left) the couple corresponding to:\\
the inter-arrival time $|$ the routing decision ($1$ means that the
customer goes back in the queue and $2$ means that the customer
leaves the network):
\begin{eqnarray*}
\begin{array}{ccccc}
1&\rightarrow&1&|&2\\
2&\rightarrow&1&|&2\\
3&\rightarrow&1&|&2\\
&\vdots&&&\\
\lfloor xn\rfloor&\rightarrow&1&|&2\\
\lfloor xn\rfloor+1&\rightarrow&1&|&1,2\\
\lfloor xn\rfloor+2&\rightarrow&1&|&2\\
&\vdots&&&\\
n&\rightarrow&1&|&2\\
n+1&\rightarrow&n&|&\underbrace{1,\dots,1}_{\lfloor
  xn\rfloor},2\\
n+2&\rightarrow&1&|&\underbrace{1,\dots,1}_{n-\lfloor
  xn\rfloor},2\\
n+3&\rightarrow&1&|&2\\
&\vdots&&&\iffalse\\

n+\lfloor xn\rfloor+1&\rightarrow&1&|&2\\
n+\lfloor xn\rfloor+2&\rightarrow&1&|&1,2\\
n+\lfloor xn\rfloor+3&\rightarrow&1&|&2\\
&\vdots&&&\\
2n+1&\rightarrow&1&|&2\\
2(n+1)&\rightarrow&n&|&\underbrace{1,\dots,1}_{\lfloor
  xn\rfloor},2\\
&\vdots&&& \fi
\end{array}
\end{eqnarray*}

\bibliographystyle{abbrv}
\bibliography{ex}

\end{document}